\documentclass[10.8pt]{amsart}
\usepackage{epsfig}
\usepackage{amsmath, amssymb, euscript,  enumerate}
\usepackage{bbm}
\usepackage{color}
\usepackage{verbatim}

\newcommand{\supp}{\text{\rm supp}}

\newcommand{\ap}{\alpha}             
\newcommand{\bt}{\beta}
\newcommand{\gm}{\gamma}             \newcommand{\Gm}{\Gamma}
\newcommand{\dt}{\delta}             
\newcommand{\vep}{\varepsilon}

\newcommand{\ld}{\lambda}            \newcommand{\Ld}{\Lambda}
\newcommand{\sm}{\sigma}             
\newcommand{\vp}{\varphi}
\newcommand{\om}{\omega}             \newcommand{\Om}{\Omega}
\newcommand{\vr}{\varrho}            \newcommand{\iy}{\infty}
            
\newcommand{\f}{\frac}             \newcommand{\el}{\ell}

\newcommand{\fF}{{\mathfrak F}}

\newcommand{\fL}{{\mathfrak L}}

\newcommand{\fm}{{\mathfrak m}}

\newcommand{\fo}{{\mathfrak o}}

\newcommand{\BI}{{\mathbb I}}
\newcommand{\BM}{{\bold M}}
\newcommand{\BN}{{\mathbb N}}

\newcommand{\BR}{{\mathbb R}}

\newcommand{\cD}{{\mathcal D}}
\newcommand{\cE}{{\mathcal E}}
\newcommand{\cF}{{\mathcal F}}

\newcommand{\cK}{{\mathcal K}}

\newcommand{\cO}{{\mathcal O}}

\newcommand{\cT}{{\mathcal T}}

\newcommand{\la}{\langle}          \newcommand{\ra}{\rangle}
\newcommand{\s}{\setminus}         
\newcommand{\n}{\nabla}            \newcommand{\e}{\eta}
\newcommand{\pa}{\partial}        \newcommand{\fd}{\fallingdotseq}
       
    \newcommand{\ds}{\displaystyle}

\newcommand{\bI}{\bold I}            \newcommand{\bM}{\bold M}
 \newcommand{\btau}{\boldsymbol\tau}
   
 \newcommand{\pf }{\noindent{\it Proof. }}
\newcommand{\rk }{\noindent{\it Remark. }}

  \newcommand{\dd }{\text{\rm d}}

\newcommand{\rB }{{\text{\rm B}}}   \newcommand{\rC }{{\text{\rm C}}}
\newcommand{\rL }{{\text{\rm L}}}

\newtheorem{thm}[subsection]{Theorem}
\newtheorem{lemma}[subsection]{Lemma}
\newtheorem{cor}[subsection]{Corollary}
\newtheorem{remark}[subsection]{Remark}

\newtheorem{definition}[subsection]{Definition}

\newtheorem{assumption}[subsection]{Assumption}

\numberwithin{equation}{section}

\title[Cordes-Nirenberg type estimates]{
Cordes-Nirenberg type estimates for nonlocal parabolic equations }
\author{ Yong-Cheol Kim and Ki-Ahm Lee }
\begin{document}
\begin{abstract} In this paper, we obtain Cordes-Nirenberg
type estimates for nonlocal parabolic equations on the more flexible
solution space $L^{\iy}_T(L^1_{\om})$ than the classical solution
space $\rB(\BR^n_T)$ consisting of all bounded functions on
$\BR^n_T$.
\end{abstract}

\address{$\bullet$ Yong-Cheol Kim : Department of Mathematics Education, Korea University, Seoul 136-701,
Republic of Korea $\,\,\,\,\,\,\&\,\,\,\,\,\,$ Department of Mathematics, Korea Institute for Advanced Study, Seoul 130-722, Republic of Korea}

\email{ychkim@korea.ac.kr}

\address{$\bullet$ Ki-Ahm Lee : Department of Mathematics, Seoul National University, Seoul 151-747,
Republic of Korea $\,\,\,\,\,\,\,\,\,\,\,\&\,\,\,\,\,\,\,\,\,\,\,$ Department of Mathematics, Korea Institute for Advanced Study, Seoul 130-722, Republic of Korea} 

\email{kiahm@math.snu.ac.kr}

\maketitle

\tableofcontents

\section{Introduction}

\subsection{Nonlocal parabolic equations.}
In this paper, we study Cordes-Nirenberg type estimates on the more
flexible solution space $L^{\iy}_T(L^1_{\om})$ for nonlocal
parabolic equations. In \cite{KL}, we obtained
interior $C^{1,\ap}$-estimates on the solution space $\rB(\BR^n_T)$
for nonlocal parabolic translation-invariant equations, and also the
reader can refer to \cite{CS1} and \cite{CS2} for the elliptic case.

Throughout this paper, we will consider the purely nonlocal parabolic
Isaacs equations of the form
\begin{equation}
\begin{split}&\bI u(x,t)-\pa_t u(x,t):=\inf_{a\in A}\sup_{b\in B}\left(L_{ab}u(x,t)-\pa_t u(x,t)\right)\\
&\quad=\inf_{a\in A}\sup_{b\in
B}\left(\int_{\BR^n}\mu_t(u,x,y)K_{ab}(x,y,t)\,dy-\pa_t u(x,t)\right)\\
&\quad=f(x,t)\,\,\text{ in
$\Om\times(-\tau,0]:=\Om_{\tau}$},\,0<\tau\le
T,\end{split}\end{equation} where $\Om$ is a bounded domain in
$\BR^n$, $\mu_t(u,x,y)=u(x+y,t)+u(x-y,t)-2u(x,t)$ and
$A,B$ are arbitrary sets. We call such $L_{ab}$ and $K_{ab}$ a {\it linear
integro-differential operator} and the kernel of the operator $L_{ab}$, respectively. Also we simply write $L$ and $K$ without indices.

We say that an operator $L$ belongs to $\fL_0$ if its corresponding
kernel $K\in\cK_0$ satisfies the uniform ellipticity assumption
\begin{equation}(2-\sm)\f{\ld}{|y|^{n+\sm}}\le K(x,y,t)\le
(2-\sm)\f{\Ld}{|y|^{n+\sm}},\,\,0<\sm<2,\,\,y\in\BR^n\s\{0\}.
\end{equation}
If $K(x,y,t)=c_{\sm}(2-\sm)|y|^{-n-\sm}$ where $c_{\sm}>0$ is the
normalization constant, then the corresponding operator is
$L=-(-\Delta)^{\sm/2}$. Also we say the operator $L\in\fL_0$
belongs to $\fL_1=\fL_1(\sm)$ if its corresponding kernel
$K\in\cK_1=\cK_1(\sm)$ satisfies $K\in C^1$ in $y$ away from the origin and
satisfies
         \begin{equation}
        \sup_{(x,t)\in\BR^n\times(-T,0]}|\n_y K(x,y,t)|\leq \frac{C_1}{|y|^{n+1+\sigma}}\,\,\text{ for any $y\in\BR^n\s\{0\}$.}         \end{equation}
Finally we say that the operator $L\in\fL_1$ belongs to
$\fL_2=\fL_2(\sm)$ if its corresponding kernel
$K\in\cK_2=\cK_2(\sm)$ satisfies $K=K(y)\in C^2$ away from the origin and
satisfies
         \begin{equation}
        |D^2 K(y)|\leq \frac{C_2}{|y|^{n+2+\sigma}}\,\,\text{ for any $y\in\BR^n\s\{0\}$.}         \end{equation}

We denote by $\om_{\sm}(y)=1/(1+|y|^{n+\sm})$ for $\sm\in(0,2)$ and
we write $\om:=\om_{\sm_0}$ for some $\sm_0\in(1,2)$, and also we
denote by $\om(Q_r)=\int_{Q_r}\om(y)\,dy$. Let $\fF$ denote the
family of all real-valued measurable functions defined on
$\BR^n_T:=\BR^n\times(-T,0]$. For $u\in\fF$ and $t\in(-T,0]$, we
define the weighted norm $\|u(\cdot,t)\|_{L^1_{\om}}$ by
$$\|u(\cdot,t)\|_{L^1_{\om}}=\int_{\BR^n}|u(x,t)|\om(x)\,dx.$$
We consider the function space $L^{\iy}_T(L^1_{\om})$ of all
continuous $L^1_{\om}$-valued functions $u\in\fF$ given by the
family
$$\biggl\{u\in\fF:\sup_{s\in(-T,0]}\|u(\cdot,s)\|_{L^1_{\om}}<\iy,
\lim_{s\to t^-}\|u(\cdot,s)-u(\cdot,t)\|_{L^1_{\om}}=0\text{ for any
$t\in(-T,0]$}\biggr\}$$ with the norm
$\|u\|_{L^{\iy}_T(L^1_{\om})}=\ds\sup_{t\in(-T,0]}\|u(\cdot,t)\|_{L^1_{\om}},$
which is separable with respect to the topology given by the norm.

A mapping $\BI:\fF\to\fF$ given by $u\mapsto\BI u$ is called a {\it
nonlocal parabolic operator} if
(a) $\BI u(x,t)$ is well-defined for any $u\in C^2_x(x,t)\cap
L^{\iy}_T(L^1_{\om})$,

\noindent(b) $\BI u$ is continuous on $\Om_{\tau}\subset\BR^n_T$, whenever
$u\in C^2_x(\Om_{\tau})\cap L^{\iy}_T(L^1_{\om})$, where
$C^2_x(x,t)$ is the class of all $u\in\fF$ whose second derivatives
$D^2 u$ in space variables exist at $(x,t)$ and $C^2_x(\Om_{\tau})$
denotes the class of all $u\in\fF$ such that $u\in C^2_x(x,t)$ for
any $(x,t)\in\Om_{\tau}$ and $\ds\sup_{(x,t)\in\Om_{\tau}}|D^2
u(x,t)|<\iy$. Such a nonlocal operator $\BI$ is said to be {\it
uniformly elliptic} with respect to a class $\fL$ of linear
integro-differential operators if
\begin{equation}\BM^-_{\fL}v(x,t)\le
\BI(u+v)(x,t)-\BI u(x,t)\le\BM^+_{\fL}v(x,t)
\end{equation} where
$\BM^-_{\fL}v(x,t):=\inf_{L\in\fL}L v(x,t)$ and
$\BM^+_{\fL}v(x,t):=\sup_{L\in\fL}L v(x,t)$.

We consider the corresponding maximal and minimal operators
\begin{equation*}\begin{split}\BM_{\fL_0}^+ u(x,t)&=\sup_{L\in\fL_0}Lu(x,t)=\int_{\BR^n}
\f{\Ld\mu_t(u,x,y)^+-\ld\mu_t(u,x,y)^-}{|y|^{n+\sm}}\,dy,\\
\BM_{\fL_0}^- u(x,t)&=\inf_{L\in\fL_0}Lu(x,t)=\int_{\BR^n}
\f{\ld\mu_t(u,x,y)^+-\Ld\mu_t(u,x,y)^-}{|y|^{n+\sm}}\,dy.
\end{split}\end{equation*}

\subsection{Outline}
In Section 2, we get various parabolic interpolation inequalities
which facilitate the required estimates for viscosity solutions for
nonlocal parabolic equations. In Section 3, we obtain H\"older
regularities and interior $C^{1,\ap}$-estimates of such viscosity
solutions by applying the result of \cite{KL} (refer to \cite{CS1}
for the elliptic case). In Section 4, we get boundary estimates and
global estimates by certain parabolic adaptation of the barrier
function which was used in \cite{CS1} for the elliptic case. In
Section 5, we establish stability properties of viscosity solutions
and it was proved that if two nonlocal parabolic equations are very
close to each other in certain sense, then so are those solutions.
The parabolic case has time shift contrary to the elliptic case, and
so this obstacle shall be overcome in this section. In Section 6, we
obtain $C^{1,\ap}$-regularity for nonlocal parabolic equations with
variable coefficients. Finally, in Section 7, we furnish a parabolic
version of the integral Cordes-Nirenberg type estimates and various
applications including $C^{2,\ap}$-regularity for nonlocal parabolic
equations.

\subsection{Notations and Definitions}
We write the notations and definitions briefly for the reader.

\begin{itemize}

\item $B_r=B_r(0)$ and $\BR^n_T=\BR^n\times(-T,0]$ for $r>0$ and
$T>0$.

\item $Q_r=B_r\times I_r^{\sm}$ and $Q_r(x,t)=Q_r+(x,t)$ for $r>0$,
$(x,t)\in\BR^n_T$ and $I_r^{\sm}=(-r^{\sm},0]$ with $\sm\in(0,2)$.

\item $\pa_p\Om_\tau:=\pa_x\Om_\tau\cup\pa_b\Om_\tau:=\pa\Om\times(-\tau,0]\cup\Om\times\{-\tau\}$
for a bounded domain $\Om\subset\BR^n$ and $\tau\in(0,T)$.

\item For $X=(x,t),Y=(y,s)\in\BR^n_T$, we define the {\it parabolic
distance} $\dd$ by $$\dd((x,t),(y,s))=\begin{cases}
(|x-y|^{\sigma}+|t-s|)^{1/\sigma}, &t\le s,\\
\iy, &t>s.
\end{cases}$$
For $X_0=(x_0,t_0)\in\BR^n_T$, we set $\rB^{\dd}(x_0,t_0)=\{(x,t)\in\BR^n_T:\dd(X,X_0)<r\}$.

\item Denote by $\rB(\Om_{\tau})$ the class of all $u\in\fF$ which is bounded in $\Om_{\tau}\subset\BR^n_T$.

\item For $a,b\in\BR$, we denote by $a\vee b=\max\{a,b\}$ and
$a\wedge b=\min\{a,b\}$.

\item We denote by $\om_n$ the surface measure of the unit sphere
$S^{n-1}$ of $\BR^n$.

\item For $(z,s)\in\BR^n_T$ and $u\in\fF$, we denote the translation
operators $\btau_z$, $\btau^s$ and $\btau_z^s$ by $\btau_z
u(x,t)=u(x+z,t)$, $\btau^s u(x,t)=u(x,t+s)$ and $\btau_z^s
u(x,t)=u(x+z,t+s)$, respectively.

\item Denote by $\e$ {\it any fixed
sufficiently small positive number}.

\item Denote by $\n u$ and $D u$ the gradient and derivatives of $u$ in space variable, respectively.

\item For two quantities $a$ and $b$, we write $a\lesssim b$ (resp.
$a\gtrsim b$) if there is a universal constant $C>0$ ({\it depending
only on $\ld,\Ld,n,\e,\sm_0$ and the constants in $(1.3)$, $(1.4)$
and $(2.8)$, but not on $\sm$}) such that $a\le C\,b$ (resp. $b\le
C\,a$).

\item For $u\in C(Q_r)$, we define
$\|u\|_{C(Q_r)}=\sup_{(x,t)\in Q_r}|u(x,t)|$. For $\ap\in(0,1]$,
$\sm\in(0,2)$ and $r>0$, we define the {\it parabolic $\ap^{th}$
H\"older seminorm} of $u$ by
$$\qquad\quad[u]_{C^{\ap}(Q_r)}=\sup_{(x,t),(y,s)\in
Q_r}\f{|u(x,t)-u(y,s)|}{(|x-y|^{\sm}+|t-s|)^{\ap/\sm}}.$$

\item Let $\BI$ be a uniformly elliptic operator in the
sense of $(1.2)$ with respect to some class $\fL$ and let
$f:\BR^n_T\to\BR$ be a continuous function. Then a function
$u\in\fF$ upper $($lower$)$ semicontinuous on $\overline\Omega\times
J$ where $J:=(a,b)\subset(-T,0]$ is said to be a {\it viscosity
subsolution $($viscosity supersolution$)$} of an equation $\BI
u-\pa_t u=f$ on $\Omega\times J$ and we write $\BI u-\pa_t u\ge f$
$($ $\BI u-\pa_t u\le f$ $)$ on $\Omega\times J$ in the viscosity
sense, if for each $(x,t)\in\Omega\times J$ there is some
neighborhood $Q_r(x,t)\subset\Omega\times J$ of $(x,t)$ such that
$\BI v(x,t)-\pa_t\vp(x,t)\ge f(x,t)$ $($ $\BI
v(x,t)-\pa_t\vp(x,t)\le f(x,t)$ $)$ for
$v=\vp\mathbbm{1}_{Q_r(x,t)}+u\mathbbm{1}_{\BR^n_T\s Q_r(x,t)}$
whenever $\vp\in C^2(Q_r(x,t))$ with $\vp(x,t)=u(x,t)$ and $\vp>u$
$($ $\vp<u$ $)$ on $Q_r(x,t)\s\{(x,t)\}$ exists. Also a function $u$
is called as a {\it viscosity solution} if it is both a viscosity
subsolution and a viscosity supersolution to $\BI u-\pa_t u=f$ on
$\Omega\times J$.

\item In what follows, we denote by $\bI_0=\inf_{L\in\fL}L$ where $\fL$ is
 a family of linear integro-differential operators in $\fL_0$.

\end{itemize}

\section{Parabolic interpolation inequalities}

Let $u\in C(Q_r)$. For $0<\ap\le 1$ and $\sm\in(0,2)$, we define the
{\it $\ap^{th}$ H\"older seminorms} of $u$ in the space and time
variable as follows, respectively;

(i) $\ds
[u]_{C^{\ap}_x(Q_r)}=\sup_{t\in(-r^{\sm},0]}\,\sup_{(x,t),(y,t)\in
Q_r}\f{|u(x,t)-u(y,t)|}{|x-y|^{\ap}},$

(ii) $\ds [u]_{C^{\ap}_t(Q_r)}=\sup_{x\in B_r}\,\sup_{(x,t),(x,s)\in
Q_r}\f{|u(x,t)-u(x,s)|}{|t-s|^{\ap}}.$

\noindent If $\,0<\ap/\sm\le 1$, then it is easy to check that the
seminorms
$[\,\cdot\,]_{C^{\ap}_x(Q_r)}+[\,\cdot\,]_{C^{\f{\ap}{\sm}}_t(Q_r)}$
and $[\,\cdot\,]_{C^{\ap}(Q_r)}$ are equivalent.

We give an useful parabolic interpolation inequality associated with
our equations.

\begin{thm} Let $\sm\in(\sm_0,2)$ with $\sm_0\in(1,2)$ and let $\fL$ be a family of linear integro-differential operators in $\fL_0$. If $u\in L^{\iy}_T(L^1_{\om})$ is a viscosity solution of
the nonlocal parabolic equation
\begin{equation*}\bI_0 u-\pa_t u=0\,\text{ in
$Q_2$,}\end{equation*} where $\bI_0$ is defined in $\fL$,
then
$\sup_{r\in(0,1)}\|u\|_{C(Q_r)}\lesssim\|u\|_{L^{\iy}_T(L^1_\om)}$.
Moreover, if $\,\bI_0$ is defined in $\fL=\fL_2(\sm)$, then any viscosity
solution  $u\in L^{\iy}_T(L^1_{\om})$ of the nonlocal parabolic
equation admits the estimate
$$\sup_{r\in(0,1)}\|(-\Delta)^{\sm/2}u\|_{C(Q_r)}\vee\,\|\pa_t u\|_{C(Q_r)}\lesssim
\|u\|_{L^{\iy}_T(L^1_\om)}.$$
\end{thm}

\pf For the first part, without loss of generality, we may assume
that $u$ is bounded on $\BR^n_T$. Indeed, if we write $u=u_1+u_2$
where $u_1=u\mathbbm{1}_{Q_r}$, then it easily follows from the
uniform ellipticity of $\bI_0$ that
$$(i)\,\,\bM^+_{\fL_0}u_1-\pa_t u_1\gtrsim-\|u\|_{L^{\iy}_T(L^1_\om)}\text{ and
}(ii)\,\,\bM^-_{\fL_0}u_1-\pa_t
u_1\lesssim\|u\|_{L^{\iy}_T(L^1_\om)}\text{ in $Q_r$.}$$

In case of (i), the estimate
$\sup_{Q_r}u\lesssim\|u\|_{L^{\iy}_T(L^1_\om)}$ can be obtained by applying a parabolic Harnack inequality \cite{KL}. Indeed, we may assume
that $\|u\|_{L^{\iy}_T(L^1_{\om})}=1$ by dividing $u$ by the norm
$\|u\|_{L^{\iy}_T(L^1_{\om})}$. Thus it suffices to show that
$\,\sup_{Q_r} u\le C.$ If $u$ is non-positive on $Q_r$, then
there is nothing to prove it. Thus we may suppose that $u$ is
non-negative on $Q_r$. We set $$s_0=\inf\{s>0:u(x,t)\le s\,
\dd((x,t),\pa_p Q_{2r})^{-n-\sm},\,\forall\,(x,t)\in Q_{2r}\}.$$ Then we
see that $s_0>0$ and there is some $(\check x,\check t)\in Q_{2r}$ such
that $$u(\check x,\check t)=s_0\,\dd((\check x,\check t),\pa_p
Q_{2r}))^{-n-\sm}=s_0\dd_0^{-n-\sm}$$ where $\dd_0=\dd((\check x,\check
t),\pa_p Q_{2r})\le (2^{\sm}+2^{\sm})^{1/\sm}r<4$ for $\sm\in(1,2)$. We note that
\begin{equation}\rB^{\dd}_\rho(x_0,t_0)\subset
Q_\rho(x_0,t_0)\subset\rB^{\dd}_{2\rho}(x_0,t_0)\end{equation} for any
$\rho>0$ and $(x_0,t_0)\in\BR^n_T$.

To finish the proof, we have only to show that $s_0$ can not be too
large because $u(x,t)\le C_1 \dd((x,t),\pa_p Q_{2r})^{-n-\sm}\le C$ for
any $(x,t)\in Q_r\subset Q_{2r}$ if $C_1>0$ is some constant with
$s_0\le C_1$. Assume that $s_0$ is very large. Then by Chebyshev's
inequality we have that
$$\bigl|\{u\ge u(\check x,\check t)/2\}\cap
Q_{2r}\}\bigr|\le\f{2}{|u(\check x,\check
t)|}\|u\|_{L^{\iy}(L^1_{\om})}\lesssim s_0^{-1}\dd_0^{n+\sm}.$$ Since
$\rB^{\dd}_\rho(\check x,\check t)\subset Q_{2r}$ and
$|\rB^{\dd}_{\rho}|=C\dd_0^{n+\sigma}$ for $\rho=\dd_0/2<2$ for $\sm\in(1,2)$, we easily obtain that
\begin{equation}\bigl|\{u\ge u(\check x,\check t)/2\}\cap
\rB^{\dd}_\rho(\check x,\check t)\}\bigr|\lesssim
s_0^{-1}|\rB^{\dd}_\rho|.\end{equation} In order to get a
contradiction, we estimate $|\{u\le u(\check x,\check t)/2\}\cap
\rB^{\dd}_{\dt \rho/2}(\check x,\check t)|$ for some very small $\dt>0$
(to be determined later). For any $(x,t)\in\rB^{\dd}_{2\dt \rho}(\check
x,\check t)$, we have that $u(x,t)\le
s_0(\dd_0-\dt\dd_0)^{-n-\sm}=u(\check x,\check t)(1-\dt)^{-n-\sm}$
for $\dt>0$ so that $(1-\dt)^{-n-\sm}$ is close to $1$. We consider
the function
$$v(x,t)=\f{u(\check x,\check t)}{(1-\dt)^{n+\sm}}-u(x,t).$$ Then we see that $v\ge 0$ on
$\rB^{\dd}_{2\dt \rho}(\check x,\check t)$, and also $\bM_0^- v-\pa_t
v\le 1$ on $Q_{\dt \rho}(\check x,\check t)$ because $\bM_0^+ u-\pa_t
u\ge -1$ on $Q_{\dt \rho}(\check x,\check t)$. In order to apply
Theorem 4.12 \cite{KL} to $v$, we consider $w=v^+$ instead of $v$.
Since $w=v+v^-$, we have that
\begin{equation}\bM_0^- w-\pa_t w\le\bM_0^- v-\pa_t v+\bM_0^+ v^--\pa_t v^-\le
1+\bM_0^+ v^--\pa_t v^-\end{equation} on $Q_{\dt \rho}(\check x,\check
t).$ Since $v^-\equiv 0$ on $\rB^{\dd}_{2\dt \rho}(\check x,\check t)$,
if $(x,t)\in Q_{\dt \rho}(\check x,\check t)$ then we have that
$\mu_t(v^-,x,y)=v^-(x+y,t)+v^-(x-y,t)$ for $y\in\BR^n$.

Take any $(x,t)\in Q_{\dt \rho}(\check x,\check t)$ and any
$\vp\in\rC^2_{Q_{\dt \rho}(\check x,\check t)}(v^-;x,t)^+$. Since
$(x,t)+Q_{\dt \rho}\subset Q_{2\dt \rho}(\check x,\check t)$ and
$v^-(x,t)=0$, we see that $\pa_t\vp(x,t)=0$. Thus we have that
\begin{equation*}\begin{split}\bM_0^+
v^-(x,t)-\pa_t\vp(x,t)&=(2-\sm)\int_{\BR^n}\f{\Ld\mu_t^+(v^-,x,y)-\ld\mu_t^-(v^-,x,y)}{|y|^{n+\sm}}\,dy\\
&\le 2(2-\sm)\Ld\int_{\{y\in\BR^n:\,v(x+y,t)<0\}}\f{-v(x+y,t)}{|y|^{n+\sm}}\,dy\\
&\le 2(2-\sm)\Ld\int_{B^c_{\dt
\rho}}\f{\bigl(u(x+y,t)-(1-\dt)^{-n-\sm}u(\check x,\check t)\bigr)_+}{|y|^{n+\sm}}\,dy\\
&\leq  C(2-\sm)\Ld\bigl((\dt
\rho)^{-n-\sigma}+1\bigr)\int_{\BR^n}\f{|u(y,t)|}{1+|y|^{n+\sm}}\,dy.
\end{split}\end{equation*}
This implies that $$\bM_0^+ v^--\pa_t v^-\lesssim
\|u\|_{L^{\iy}_T(L^1_{\om})}(\dt \rho)^{-n-\sigma}\lesssim(\dt
\rho)^{-n-\sigma}\text{ on $Q_{\dt \rho}(\check x,\check t)$. }$$ Thus by
(2.3), we obtain that $w$ satisfies
$$\bM_0^- w(x,t)-\pa_t w\lesssim (\dt
\rho)^{-n-\sm}\,\,\text{ on $Q_{\dt \rho}(\check x,\check t)$ }$$ in
viscosity sense. Since $u(\check x,\check
t)=s_0\dd_0^{-\bt}=2^{-\bt}s_0 \rho^{-\bt}$, by Theorem 4.12 \cite{KL}
there is some $\vep_*>0$ such that
\begin{equation*}\begin{split}&\bigl|\{u\le u(\check
x,\check t)/2\}\cap\rB^{\dd}_{\dt \rho/2}(\check x,\check
t)\bigr|\le\bigl|\{u\le u(\check x,\check t)/2\}\cap Q_{\dt
\rho/2}(\check x,\check t)\bigr|\\
&\qquad\qquad\qquad=\bigl|\{w\ge u(\check x,\check
t)((1-\dt)^{-\bt}-1/2)\}
\cap Q_{\dt \rho/2}(\check x,\check t)\bigr|\\
&\qquad\qquad\qquad\lesssim(\dt
\rho)^{n+\sigma}\bigl[((1-\dt)^{-\bt}-1)u(\check x,\check t)+C(\dt
\rho)^{-\sm}(\dt
\rho)^{\sm}\bigr]^{\vep_*}\\
&\qquad\qquad\qquad\qquad\qquad\qquad\qquad\qquad\times\bigl[u(\check x,\check t)((1-\dt)^{-\bt}-1/2)\bigr]^{-\vep_*}\\
&\qquad\qquad\qquad\lesssim(\dt
\rho)^{n+\sigma}\biggl[\biggl(\f{(1-\dt)^{-\bt}-1}{(1-\dt)^{-\bt}-1/2}\biggr)^{\vep_*}
+\f{s_0^{-\vep_*}\rho^{n+\sm}}{((1-\dt)^{-\bt}-1/2)^{\vep_*}}\biggr]\\
&\qquad\qquad\qquad\lesssim(\dt
\rho)^{n+\sigma}[((1-\dt)^{-\bt}-1)^{\vep_*}+s_0^{-\vep_*}\rho^{n+\sm}].
\end{split}\end{equation*}
We now choose $\dt>0$ so small enough that $(\dt
\rho)^{n+\sigma}((1-\dt)^{-\bt}-1)^{\vep_*}\lesssim |\rB^{\dd}_{\dt
\rho/2}|/4.$ Since $\dt$ was chosen independently of $s_0$, if $s_0$ is
large enough for such fixed $\dt$ then we get that $(\dt
\rho)^{n+\sm}s_0^{-\vep_*}\rho^{n+\sm}\lesssim |\rB^{\dd}_{\dt \rho/2}|/4.$
Therefore we obtain that $$\bigl|\{u\le u(\check x,\check
t)/2\}\cap\rB^{\dd}_{\dt \rho/2}(\check x,\check t)\bigr|\le
|\rB^{\dd}_{\dt \rho/2}|/2.$$ Thus we conclude that
\begin{equation*}\begin{split}\bigl|\{u\ge u(\check x,\check t)/2\}\cap\rB^{\dd}_\rho(\check x,\check t)\bigr|
&\ge\bigl|\{u\ge u(\check x,\check t)/2\}\cap \rB^{\dd}_{\dt
\rho/2}(\check x,\check t)\bigr|\\&\ge\bigl|\{u>u(\check x,\check
t)/2\}\cap\rB^{\dd}_{\dt \rho/2}(\check x,\check
t)\bigr|\\&\ge\bigl|\rB^{\dd}_{\dt \rho/2}(\check x,\check
t)\bigr|-\bigl|\rB^{\dd}_{\dt \rho/2}\bigr|/2\\&=\bigl|\rB^{\dd}_{\dt
\rho/2}\bigr|/2=C |\rB^{\dd}_\rho|,\end{split}\end{equation*} which
contradicts (2.2) if $s_0$ is large enough.

Since $-u$ is another solution of the given equation, the second
equation (ii) can be transformed equivalently to the equation
$$\bM^+_{\fL_0}(-u_1)-\pa_t(-u_1)\ge-\|u\|_{L^{\iy}_T(L^1_\om)}\text{ in
$Q_r$.}$$ This and the upper bound estimate in the above imply that
$\inf_{Q_r}u\gtrsim-\,\|u\|_{L^{\iy}_T(L^1_\om)}$.

The second part can be obtained from Corollary 7.3 in \cite{KL1}.
Hence we complete the proof. \qed

Next we prove various lemmas which furnish parabolic interpolation
inequalities.

\begin{lemma} If $\,u\in L^{\iy}_T(L^1_{\om})$ is a function with
$u(\cdot,t)\in C^k(B_r)$ for $t\in(-r^{\sm},0]$ and
$[D^{\bt}u]_{C_x^{\ap}(Q_r)}<\iy$ for some $\ap\in(0,1)$, then for
each $t\in(-r^{\sm},0]$ and multiindex $\bt$ with $|\bt|=k\in\BN$,
there exists some $z^t_0\in B_r$ $($depending on $t$$)$ such that
$$\bigl|D^{\bt}u(z^t_0,t)\bigr|\le
\bigl(\f{3r}{2}\bigr)^{\ap}\,[D^{\bt}u]_{C_x^{\ap}(Q_r)}+\f{2(4k)^k}{\om(B_{r/2})\,r^k}\|u\|_{L^{\iy}_T(L^1_{\om})}.$$
\end{lemma}

\pf Take $h=\f{r}{2k}$ and any multiindex $\bt$ with $|\bt|=k$. For
$(y,t)\in B_{r/2}\times(-T,0]$, we consider the finite difference
operator $\,D^{\bt}_h u(y,t)=D^{\bt_1}_{h,1}\,D^{\bt_2}_{h,2}\cdots
D^{\bt_n}_{h,n}u(y,t)$ where
$$D_{h,i}u(y,t)=\f{1}{h}\,[u(y+h e_i,t)-u(y,t)]$$
for a standard basis $\{e_1,\cdots,e_n\}$ of $\BR^n$. For
$i=1,\cdots,n$, we observe that
\begin{equation}D^{\bt_i}_{h,i}u(y,t)=\f{1}{h^{\bt_i}}\sum_{s=0}^{\bt_i}(-1)^{s}\f{\bt_i
!}{(\bt_i-s)!\,s !}\,u\bigl(y+(\bt_i-s)h e_i,t\bigr).
\end{equation}
By the mean value theorem, we see that there are some $z^t_1\in
B_h(y)$ and $z^t_2\in B_{2h}(y)$ such that
$$D_{h,i}D_{h,j}u(y,t)=\pa_{y_i}[D_{h,j}u](z^t_1,t)=D_{h,j}(\pa_{y_i}u)(z^t_1,t)
=\pa_{y_i y_j}u(z^t_2,t).$$ This implies that $D^{\bt}_h
u(y,t)=D^{\bt}u(z^t_y,t)$ for some $z^t_y\in B_{r/2}(y)$. Thus it
follows from this and (2.4) that
\begin{equation*}\begin{split}
\om(B_{r/2})\bigl|D^{\bt}u(z^t_0,t)\bigr|&\le\biggl|\om(B_{r/2})D^{\bt}u(z^t_0,t)-\int_{\BR^n}D_h^{\bt}u(y,t)\,\om(y)\,dy\biggr|
+\f{2^k}{h^k}\|u\|_{L^{\iy}_T(L^1_{\om})}\\
&\le\int_{B_{r/2}}\bigl|D^{\bt}u(z^t_0,t)-D^{\bt}u(z^t_y,t)\bigr|\,\om(y)\,dy+\f{2^{k+1}}{h^k}\|u\|_{L^{\iy}_T(L^1_{\om})}\\
&\le[D^{\bt}u]_{C_x^{\ap}(Q_r)}
\bigl(\f{3r}{2}\bigr)^{\ap}\,\om(B_{r/2})+\f{2(4k)^k}{r^k}\|u\|_{L^{\iy}_T(L^1_{\om})}.
\end{split}\end{equation*}
Therefore, this completes the proof. \qed

In order to understand the {\bf parabolic H\"older spaces}
$C^{k,\gm}(Q_r)$ with $k\in\BN$ and $\gm\in(0,1)$, we define the
H\"older spaces $C^{k,\gm}_x(Q_r)$ and $C^{k,\gm}_t(Q_r)$ in the
space and time variable, respectively. For $u\in C(Q_r)$, we define
the norms
\begin{equation*}\begin{split}\|u\|_{C^{k,\gm}_x(Q_r)}&=\|u\|_{C(Q_r)}
+\sum_{i=1}^k\|D^i u\|_{C(Q_r)}+[D^k u]_{C^{\gm}_x(Q_r)},\\
\|u\|_{C^{k,\gm}_t(Q_r)}&=\|u\|_{C(Q_r)}+\sum_{i=1}^k\|\pa_t^i
u\|_{C(Q_r)} +[\pa_t^k
u]_{C^{\gm}_t(Q_r)},\end{split}\end{equation*} where $\|D^i
u\|_{C(Q_r)}=\sum_{|\bt|=i}\|D^{\bt}u\|_{C(Q_r)}$ and $[D^k
u]_{C^{\gm}_x(Q_r)}=\sum_{|\bt|=k}[D^{\bt}u]_{C^{\gm}_x(Q_r)}$ for
$i,k\in\BN$. And we denote by
$C_x^{k,\gm}(Q_r)=\{u\in\fF(\BR^n_T):\|u\|_{C^{k,\gm}_x(Q_r)}<\iy\}$
and
$C_t^{k,\gm}(Q_r)=\{u\in\fF(\BR^n_T):\|u\|_{C^{k,\gm}_t(Q_r)}<\iy\}$.

\begin{thm} If $\,u\in L^{\iy}_T(L^1_{\om})$ is a function such that
$u(\cdot,t)\in C^k(B_r)$ for each $t\in(-r^{\sm},0]$ and
$\,\sup_{|\bt|=k}[D^{\bt}u]_{C_x^{\ap}(Q_r)}<\iy$ for some
$\ap\in(0,1)$, then we have that
$$\|D^k u\|_{C(Q_r)}\le
2\,\bigl(\f{3r}{2}\bigr)^{\ap}\,[D^k u]_{C_x^{\ap}(Q_r)}+\f{2\,
c_k(4k)^k}{\om(B_{r/2})\,r^k}\|u\|_{L^{\iy}_T(L^1_{\om})},\text{
where $c_k=\sum_{|\bt|=k}1$.}$$
\end{thm}

\pf From Lemma 2.2, for any $(x,t)\in Q_r$ we obtain that
\begin{equation*}\begin{split}\bigl|D^{\bt}u(x,t)\bigr|&\le\bigl|D^{\bt}u(x,t)-D^{\bt}u(z^t_0,t)\bigr|
+\bigl|D^{\bt}u(z^t_0,t)\bigr|\\
&\le 2\,[D^{\bt}u]_{C_x^{\ap}(Q_r)}
\bigl(\f{3r}{2}\bigr)^{\ap}+\f{2(4k)^k}{\om(B_{r/2})\,r^k}\|u\|_{L^{\iy}_T(L^1_{\om})}.
\end{split}\end{equation*}
Taking the supremum on $Q_r$ and adding up on the multiindices $\bt$
with $|\bt|=k$ in the above, we easily obtain the required result.
\qed

If $\sm\in(\sm_0,2)$ for $\sm_0\in(1,2)$ and $\ap\in(0,\sm_0-1)$,
then $0<\ap<2+\ap-\sm<1$ and
\begin{equation}\f{2+\ap-\sm}{\sm}+1=\f{2+\ap}{\sm}.
\end{equation}
Then we define the {\bf parabolic H\"older space} $C^{2,\ap}(Q_r)$
endowed with the norm
\begin{equation*}\begin{split}\|u\|_{C^{2,\ap}(Q_r)}&=\|u\|_{C(Q_r)}
+\sum_{i=1}^2\|D^i u\|_{C(Q_r)}+\|\pa_t u\|_{C(Q_r)}\\
&\quad+[D^2 u]_{C^{\ap}(Q_r)}+[\pa_t u]_{C^{2+\ap-\sm}(Q_r)}.
\end{split}\end{equation*}

In the same case as the above, we can learn from Theorem 2.1 and
Theorem 2.3 that the estimates on the norm $\|u\|_{C^{2,\ap}(Q_r)}$
must be controlled by those on the seminorms $[\pa_t
u]_{C^{2+\ap-\sm}(Q_r)}\sim [\pa_t u]_{C_x^{2+\ap-\sm}(Q_r)}+[\pa_t
u]_{C_t^{\f{2+\ap-\sm}{\sm}}(Q_r)}$ and $\ds[D^2
u]_{C^{\ap}(Q_r)}\sim[D^2 u]_{C_x^{\ap}(Q_r)} +[D^2
u]_{C_t^{\f{\ap}{\sm}}(Q_r)}$. Similarly, the other parabolic
H\"older spaces can be defined along this line.

\begin{lemma} Let $\sm\in[\sm_0,2)$ for some $\sm_0\in(1,2)$ and $\ap\in(0,\sm_0-1)$ as in $(2.2)$.
If $u\in L^{\iy}_T(L^1_{\om})$ is a function with $u(x,\cdot)\in
C^1(-r^{\sm},0]$ for $x\in B_r$ and $[\pa_t
u]_{C_t^{\f{2+\ap-\sm}{\sm}}(Q_r)}<\iy$, then we have that
\begin{equation*}\|\pa_t u\|_{C(Q_r)}\le r^{2+\ap-\sm}[\pa_t
u]_{C_t^{\f{2+\ap-\sm}{\sm}}(Q_r)}+\f{4}{r^{\sm}}\|u\|_{C(Q_r)}.
\end{equation*}
\end{lemma}

\pf Take any $r\in(0,2)$ and $(x,t)\in Q_r$. Then there is some
$t_0\in(-r^{\sm},0]$ such that $|t-t_0|=\f{1}{2}r^{\sm}$, and by the
mean value theorem, there is some $t_0^x$ between $t$ and $t_0$ such
that $u(x,t_0)-u(x,t)=\f{1}{2}r^{\sm}\,\pa_t u(x,t_0^x)$. Thus we
have the estimate
\begin{equation*}\begin{split}
\f{1}{2}\,r^{\sm}\,|\pa_t u(x,t)|&\le\biggl|\f{1}{2}\,r^{\sm}\,\pa_t
u(x,t)-\bigl(u(x,t_0)-u(x,t)\bigr)\biggr|+2\,\|u\|_{C(Q_r)}\\
&=\f{1}{2}\,r^{\sm}\bigl|\pa_t u(x,t)-\pa_t
u(x,t_0^x)\bigr|+2\,\|u\|_{C(Q_r)}\\
&\le\f{1}{2}\,r^{2+\ap}[\pa_t
u]_{C_t^{\f{2+\ap-\sm}{\sm}}(Q_r)}+2\,\|u\|_{C(Q_r)}.
\end{split}\end{equation*}
Hence this implies the required inequality. \qed

\begin{lemma} Let $\sm\in[\sm_0,2)$ for some $\sm_0\in(1,2)$, and
let $u\in L^{\iy}_T(L^1_{\om})$ be a viscosity solution of the
equation
$$\bI_0 u-\pa_t u=0\,\,\text{ in $Q_2$}$$ where $\bI_0$ is defined on $\fL_2(\sm)$. If $u\in C^{2,\ap}(Q_2)$,
then we have the estimates
\begin{equation*}\begin{split}[D^2 u]_{C_t^{\f{\ap}{\sm}}(Q_r)}
&\lesssim\|D^2 u\|_{C(Q_r)}+\|u\|_{L^{\iy}_T(L^1_{\om})},\\
[\pa_t u]_{C_x^{2+\ap-\sm}(Q_r)}&\lesssim
\|u\|_{L^{\iy}_T(L^1_{\om})}
\end{split}\end{equation*} for any $r\in(0,1)$.
\end{lemma}

\pf Take any $r\in(0,1)$ and $(x,t)\in Q_r$. We consider the
difference quotients in the $x$-direction
$$u^h(x,t)=\f{u(x+h,t)-u(x,t)}{|h|}.$$ Write $u^h=u^h_1+u^h_2$ where
$u^h_1=u^h\mathbbm{1}_{Q_r}$. By Theorem 2.4 \cite{KL}, we have that
$\bM^+_{\fL_2}u^h-\pa_t u^h\ge 0$ and $\bM^-_{\fL_2}u^h-\pa_t u^h\le
0$ on $Q_r$. Since $\pa_t u^h_2\equiv 0$ in $Q_r$, it follows from
the uniform ellipticity of $\bI_0$
with respect to $\fL_2$ that
$$\bM^+_{\fL_0}u^h_1-\pa_t u^h_1\ge -\bM^+_{\fL_2}u^h_2\,\,\text{ and
}\,\,\bM^-_{\fL_0}u^h_1-\pa_t u^h_1\le -\bM^-_{\fL_2}u^h_2\,\,\text{
in $Q_r$.}$$ Then it is easy to show that
$|\bM^+_{\fL_2}u^h_2|\vee|\bM^-_{\fL_2}u^h_2|\lesssim
\|u\|_{L^{\iy}_T(L^1_{\om})}$ in $Q_r$. So we have that
$$\bM^+_{\fL_0}u^h_1-\pa_t u^h_1\gtrsim
-\|u\|_{L^{\iy}_T(L^1_{\om})}\,\,\text{ and
}\,\,\bM^-_{\fL_0}u^h_1-\pa_t u^h_1\lesssim
\|u\|_{L^{\iy}_T(L^1_{\om})}\,\,\text{ in $Q_r$.}$$ We now consider
another difference quotients in the $x$-direction
$$w^h(x,t)=\f{u^h_1(x+h,t)-u^h_1(x,t)}{|h|}.$$
Applying Theorem 2.4 \cite{KL} again, we obtain that
$$\bM^+_{\fL_0}w^h-\pa_t w^h\gtrsim
-\|u\|_{L^{\iy}_T(L^1_{\om})}\,\text{ and
}\,\,\bM^-_{\fL_0}w^h-\pa_t w^h\lesssim
\|u\|_{L^{\iy}_T(L^1_{\om})}\,\text{ in $Q_r$.}$$ From the H\"older
estimate(Theorem 3.4) below, we get the estimate
$$[w^h]_{C^{\f{\ap}{\sm}}_t(Q_r)}\le[w^h]_{C^{\ap}(Q_r)}\lesssim
\|w^h\|_{C(Q_r)}+\|w^h\|_{L^{\iy}_T(L^1_{\om})}+\|u\|_{L^{\iy}_T(L^1_{\om})}.$$
By the mean value theorem, we easily have that
$\|w^h\|_{C(Q_r)}\le\|D^2 u\|_{C(Q_r)}$. Since
$|D\om(y,s)|+|D^2\om(y,s)|\le c\,\om(y)$, it follows from the
integration by parts that
\begin{equation}\|w^h\|_{L^{\iy}_T(L^1_{\om})}\lesssim\|u\|_{L^{\iy}_T(L^1_{\om})}.\end{equation}
Thus we obtain that
$$[w^h]_{C^{\f{\ap}{\sm}}_t(Q_r)}\lesssim
\|D^2 u\|_{C(Q_r)}+\|u\|_{L^{\iy}_T(L^1_{\om})}.$$ Taking the limit
$|h|\to 0$ in the above, we conclude that the first inequality
holds.

Take any $(x,t)\in Q_r$. Then by the uniform ellipticity we have that
\begin{equation}\begin{split}\bM^-_2(\btau^t_x u-\btau^t u)(0,0)&\le\pa_t u(x,t)-\pa_t u(0,t)\\&=\bI u(x,t)-\bI u(0,t)
\le\bM^+_2(\btau^t_x u-\btau^t u)(0,0)
\end{split}\end{equation}
Let $\vp\in C^{\iy}_c(\BR^n)$ satisfy that $\vp=1$ in
$B_1$, $\vp=0$ in $\BR^n\s B_{3/2}$ and $0\le\vp\le 1$ in $\BR^n$, and take any $\rL\in\fL_2$. Then it follows from the change of variable, the mean value theorem and (1.3) that
\begin{equation}\begin{split}\rL(\btau_x^t u-\btau^t u)(0,0)&=\int_{\BR^n}\bigl[\mu_t(u,x,y)-\mu_t(u,0,y)\bigr]\vp(y)K(y)\,dy\\
&+\int_{\BR^n}\bigl[\mu_t(u,x,y)-\mu_t(u,0,y)\bigr](1-\vp(y))K(y)\,dy\\
&\lesssim\vp^+ u(x,0)+\|u\|_{L^{\iy}_T(L^1_{\om})}\,|x|
\end{split}\end{equation}
where $$\vp^+ u(x,0)=\sup_{t\in(-T,0]}\sup_{K\in\cK_2}\int_{\BR^n}\bigl[\mu_t(u,x,y)-\mu_t(u,0,y)\bigr]\vp(y)K(y)\,dy.
$$
Similarly we can obtain that 
\begin{equation}\begin{split}\rL(\btau_x^t u-\btau^t u)(0,0)\gtrsim\vp^- u(x,0)-\|u\|_{L^{\iy}_T(L^1_{\om})}\,|x|
\end{split}\end{equation}
where $$\vp^- u(x,0)=\inf_{t\in(-T,0]}\inf_{K\in\cK_2}\int_{\BR^n}\bigl[\mu_t(u,x,y)-\mu_t(u,0,y)\bigr]\vp(y)K(y)\,dy.
$$ The estimates (2.7), (2.8) and (2.9) imply that
\begin{equation}\begin{split}\vp^- u(x,0)-\|u\|_{L^{\iy}_T(L^1_{\om})}\,|x|&\lesssim\bM^-_2(\btau_x^t u-\btau^t u)(0,0)\\&\le\pa_t u(x,t)-\pa_t u(0,t)\\
&\le\bM^+_2(\btau_x^t u-\btau^t u)(0,0)\\
&\lesssim\vp^+ u(x,0)+\|u\|_{L^{\iy}_T(L^1_{\om})}\,|x|.
\end{split}\end{equation}
Applying the method in Lemma 9.2 \cite{CS3}, we have that
$$|\vp^- u(x,0)|\vee|\vp^+ u(x,0)|\lesssim\|u\|_{L^{\iy}_T(L^1_{\om})}\,|x|^{\bt}$$
for some $\bt\in(0,1)$. Here, without loss of generality, we may assume that $\bt=2+\ap-\sm$ by applying a standard telescopic argument \cite{CC}.
Hence the second inequality can be achieved from a standard translation argument. Therefore we complete the proof. \qed

{\bf Remark 2.1.} We learned from the interpolation results obtained
in this section that the norm $\|u\|_{C^{2,\ap}(Q_r)}$ of viscosity
solutions $u\in L^{\iy}_T(L^1_{\om})$ of the equation $$\bI_0 u-\pa_t u=0\,\,\text{ in $Q_2$}$$ where $\bI_0$ is defined on $\fL_2(\sm)$ for $\sm\in[\sm_0,2)$ with some $\sm_0\in(1,2)$ is controlled by only two
seminorms $[\pa_t u]_{C_t^{\f{2+\ap-\sm}{\sm}}(Q_r)}$ and $[D^2
u]_{C_x^{\ap}(Q_r)}$, and so only two norms
$\|u\|_{C^{1,\f{2+\ap-\sm}{\sm}}_t(Q_r)}$ and
$\|u\|_{C^{2,\ap}_x(Q_r)}$.

Finally, we are going to define another {\bf parabolic H\"older
space} $C^{1,\ap}(Q_r)$ in case that $1<\sm<2$. From \cite{KL}, such
$\ap>0$ could be chosen so that $\ap<\sm-1$, i.e.
$0<\theta=\theta(\sm,\ap)=\f{1+\ap}{\sm}<1$. We learned from the
definition of $C^{2,\ap}(Q_r)$ that one derivative in time variable
amounts to two derivatives in space variable. Since there is only
one derivative in space variable on $C^{1,\ap}(Q_r)$, the space
should be defined as the family of all functions $u\in\fF$ with the
norm
$$\|u\|_{C^{1,\ap}(Q_r)}=\|u\|_{C(Q_r)}+\|D u\|_{C(Q_r)}+[D
u]_{C^{\ap}(Q_r)}+[u]_{C_t^{\theta}(Q_r)}<\iy.$$

We define the class $\fL_*$ of operators $L$ with kernels
$K\in\cK_*$ satisfying (1.2) such that there are some $\vr_0>0$ and
a constant $C>0$ such that
\begin{equation}\sup_{(x,t)\in\BR^n_T}|\n_y K(x,y,t)|\le
C\,\om(y)\,\,\text{ for any $y\in\BR^n\s B_{\vr_0}.$ }\end{equation}
We note that $\fL_1$ is the largest scale invariant class contained
in the class $\fL_*$.

\begin{thm} Let $\sm\in[\sm_0,2)$ for some $\sm_0\in(1,2)$.
Then there is some $\vr_0>0$ $($depending on $\ld,\Ld,\sm_0$ and
$n$$)$ so that if $\BI$ is a nonlocal, translation-invariant and
uniformly elliptic operator with respect to $\fL_*$ and $u\in
L^{\iy}_T(L^1_{\om})$ is a viscosity solution of the equation 
$$\BI u-\pa_t u=0\,\,\text{ in $Q_2$,}$$ then
there is some $\ap\in(0,1)$ such that
$$\|D u\|_{C_t^{\f{\ap}{\sm}}(Q_r)}\lesssim
\|D u\|_{C(Q_r)}+\|u\|_{L^{\iy}_T(L^1_{\om})}$$ for any
$r\in(0,1)$.\end{thm}

\pf We proceed the proof by applying Theorem 3.4 below to the
difference quotients in the $x$-direction
$$w^h(x,t)=\f{u(x+h,t)-u(x,t)}{|h|}.$$
Take any $r\in(0,2)$. Then we write $w^h=w_1^h+w_2^h$ where
$w_1^h=w^h\mathbbm{1}_{Q_r}$. From Theorem 2.4 \cite{KL}, we have
that $\bM^+_{\fL^*}w^h-\pa_t w^h\ge 0$ and $\bM^-_{\fL^*}w^h-\pa_t
w^h\le 0$ in $Q_r$. Because $\pa_t w_2^h\equiv 0$ in $Q_r$, it
follows from the uniform ellipticity of $\bI$ with respect to $\fL^*$ that
\begin{equation*}\begin{split} \BM^+_{\fL_0}w_1^h-\pa_t w_1^h&\ge
\BM^+_{\fL_*}w_1^h-\pa_t w_1^h\ge\BM^+_{\fL_*} w^h-\BM^+_{\fL_*}
w_2^h-\pa_t w^h\ge-\BM^+_{\fL_*}w_2^h\,\,\text{ in $Q_r$},\\
\BM^-_{\fL_0}w_1^h-\pa_t w_1^h&\le\BM^-_{\fL_*}w_1^h-\pa_t
w_1^h\le\BM^-_{\fL_*} w^h-\BM^-_{\fL_*} w_2^h-\pa_t
w^h\le-\BM^-_{\fL_*}w_2^h\,\,\text{ in
$Q_r$}.\end{split}\end{equation*} If we can show that
$|\BM^+_{\fL_*}w_2^h|\vee|\BM^-_{\fL_*}w_2^h|\lesssim
\|u\|_{L^{\iy}_T(L^1_{\om})}$ in $Q_r$, then we have that
$$\BM^+_{\fL_0}w_1^h-\pa_t w_1^h\gtrsim-\|u\|_{L^{\iy}_T(L^1_{\om})}\text{ and
}\BM^-_{\fL_0}w_1^h-\pa_t w_1^h\lesssim
\|u\|_{L^{\iy}_T(L^1_{\om})}\text{ in $Q_r$}$$ for $h$ with a
sufficiently small $|h|$. Indeed, by using (2.11), it can be obtained
from the fact that
\begin{equation*}\begin{split}&\int_{\BR^n\s
B_{\rho}}|u(x+y,t)|\f{|K(x,y,t)-K(x,y-h,t)|}{|h|}\,dy\\
&\qquad+\int_{\BR^n\s B_{\rho}}|u(x+y+h,t)|K(x,y,t)\,dy\lesssim
\|u\|_{L^{\iy}_T(L^1_{\om})}\end{split}\end{equation*} for some
$\rho>0$. Hence $w_1^h$ admits the H\"older estimate(Theorem 3.4)
below on $Q_r$, and thus applying the mean value theorem and
integration by parts with (2.11) gives the estimate
\begin{equation*}\begin{split}\|w_1^h\|_{C_t^{\f{\ap}{\sm}}(Q_r)}
&\le\|D u\|_{C(Q_r)}+\|u\|_{L^{\iy}_T(L^1_{\om})}.
\end{split}\end{equation*}
Finally, taking the limit $|h|\to 0$ in the above, we obtain the
required result. \qed

{\bf Remark 2.2.} From Theorem 2.6, we saw that the norm
$\|u\|_{C^{1,\ap}(Q_r)}$ of viscosity solutions $u\in
L^{\iy}_T(L^1_{\om})$ of the equation  $$\BI u-\pa_t u=0\,\,\text{ in $Q_2$ }$$ where  $\BI$ is a nonlocal, translation-invariant and uniformly elliptic with respect to
$\fL_*$ is controlled by only two seminorms $[D u]_{C_x^{\ap}(Q_r)}$
and $[u]_{C_t^{\theta}(Q_r)}$ with $\theta=\f{1+\ap}{\sm}$. Thus the
norm  $\|u\|_{C^{1,\ap}(Q_r)}$ is completely governed by only two
norms $\|u\|_{C_x^{1,\ap}(Q_r)}$ and $\|u\|_{C_t^{\theta}(Q_r)}$.

\section{Preliminary estimates}

In this paper, we always impose the following assumptions on $\om$;
\begin{equation}1+|y|\in L^1_{\om},\end{equation}
\begin{equation}\sup_{B_r(y)}\om\le C_r\om(y).\end{equation}
The uniform ellipticity (1.2) depends on a class $\fL$ of {\it
linear integro-differential operators}. Such an operator $L$ in
$\fL$ is of the form $L u(x,t)=\int_{\BR^n}\mu_t(u,x,y)K(x,y,t)\,dy$
for a nonnegative symmetric kernel $K$ satisfying
$$\sup_{(x,t)\in\BR^n_T}\int_{\BR^n}(1\wedge|y|^2)K(x,y,t)\,dy\le C<\iy.$$ Here
the symmetric property means that for each $(x,t)\in\BR^n_T$,
$K(x,-y,t)=K(x,y,t)$ for all $y\in\BR^n$.

Let $\Om$ be a bounded domain in $\BR^n$. Then we say that a
function $u:\BR^n_T\to\BR$ is Lipschitz in space on
$\Om_{\tau}=\Om\times(-\tau,0]$, $\tau\in(0,T)$ ( we write $u\in
C^{0,1}_x(\Om_{\tau})$ ), if there is some constant $C>0$
( independent of $x,y$ ) such that
\begin{equation}\sup_{t\in(-\tau,0]}|u(x,t)-u(y,t)|\le C|x-y|
\end{equation} for any $x,y\in\Om$. We denote by
$[u]_{C^{0,1}_x(\Om_{\tau})}$ the smallest $C$ satisfying (3.3).

We say that a function $u:\BR^n_T\to\BR$ is in
$C^{1,1}_x(\Om_{\tau})$, if there is a constant $C_0>0$ ( independent
of $(x,t)$ and $(y,t)$ ) such that
\begin{equation}|u(y,t)-u(x,t)-(y-x)\cdot\n u(x,t)|\le C_0|y-x|^2\end{equation} for all
$(x,t),(y,t)\in\Om_{\tau}$. We denote by the norm
$\|u\|_{C^{1,1}_x(\Om_{\tau})}$ the smallest $C_0$ satisfying (3.4).

The following definition is the parabolic setting of that \cite{CS1}
of the elliptic case.

\begin{definition} For a nonlocal parabolic operator $\BI$ and $\tau\in(0,T]$, we define
$\|\BI\|$ in $\Om_{\tau}$ with respect to some weight $\om$ as
$$\|\BI\|=\sup_{(y,s)\in\Om_{\tau}}\sup_{u\in\cF^{M}_{y,s}}
\f{|\BI
u(y,s)|}{1+\|u\|_{L^{\iy}_T(L^1_{\om})}+\|u\|_{C^{1,1}_x(Q_1(y,s))}}$$
where $\cF^{M}_{y,s}=\{u\in\fF\cap
C^2_x(y,s):\|u\|_{L^{\iy}_T(L^1_{\om})}+\|u\|_{C^{1,1}_x(Q_1(y,s))}\le
M\}$ for some $M>0$.
\end{definition}

\begin{assumption} If $K_{\fL}:=\sup_{\ap}K_{\ap}$ is the supremum
of all kernels corresponding to operators in the class $\fL$, then
for each $r>0$ there is a constant $C_r>0$ such that
$\sup_{(x,t)\in\BR^n_T}K_{\fL}(x,y,t)\le C_r\,\om(y)$ for any
$y\in\BR^n\s B_r$.
\end{assumption}

\begin{assumption} There is some $C>0$ such that $\,\,\sup_{L\in\fL}\|L\|\le C<\iy$.
\end{assumption}

\rk Assumption 3.3 implies that $\|\BM^+_{\fL}\|\le C$ and
$\|\BM^-_{\fL}\|\le C$.

\begin{thm} Let $\sm\in(\sm_0,2)$ for some
$\sm_0\in(1,2)$. If $u\in L^{\iy}_T(L^1_{\om})$ is a function
satisfying
$$\BM^+_{\fL_0}u-\pa_t u\ge-C_0\,\,\text{ and }\,\,\BM^-_{\fL_0}u-\pa_t u\le
C_0\,\,\,\text{ in $Q_{1+\e}$, }$$ then there is some $\ap>0$ such that
$$\|u\|_{C^{\ap}(Q_1)}\lesssim
\|u\|_{C(Q_{1+\e})}+\|u\|_{L^{\iy}_T(L^1_{\om})}+C_0.$$
\end{thm}

\pf We note that $u$ is continuous on $\overline Q_{1+\e}$. Set
$v=u\mathbbm{1}_{Q_{1+\e}}$ and $w=u\mathbbm{1}_{\BR^n_T\s
Q_{1+\e}}$. Since $u=v+w$ and $\pa_t w\equiv 0$ on $Q_1$, we have
that
\begin{equation*}\begin{split}(\BM^+_{\fL_0}v-\pa_t v)+\BM^+_{\fL_0}w&\ge
\BM^+_{\fL_0}u-\pa_t u\ge-C_0\,\text{ in $Q_1$,}\\
(\BM^-_{\fL_0}v-\pa_t v)+\BM^-_{\fL_0}w&\le\BM^-_{\fL_0}u-\pa_t u\le
C_0\,\,\,\,\,\,\text{ in $Q_1$}.
\end{split}\end{equation*}
So it suffices to show that if $(x,t)\in Q_1$, then $|L w(x,t)|\lesssim
\|u\|_{L^{\iy}_T(L^1_{\om})}$ for any $L\in\fL_0$, i.e. we have
only to show that if $(x,t)\in Q_1$, then
$$\biggl|\int_{|y|\ge 1+\e}u(y,t)K(x,x\pm y,t)\,dy\biggr|\lesssim
\|u\|_{L^{\iy}_T(L^1_{\om})}$$ for any $L\in\fL_0$. Indeed, we
note that $|y|>(1+\e)|x|$ for any $x\in B_1$ and $y\in\BR^n\s
B_{1+\e}$, and so $|x\pm y|\ge|y|-|x|\ge\f{\e}{1+\e}|y|$. Thus we
have the estimate
$$\biggl|\int_{|y|\ge 1+\e}u(y,t)K(x,x\pm
y,t)\,dy\biggr|\lesssim\int_{|y|\ge 1+\e}\f{|u(y,t)|}{|x\pm
y|^{n+\sm}}\,dy\lesssim\|u\|_{L^{\iy}_T(L^1_{\om})}.$$ This implies
that
$$\BM^+_{\fL_0}v-\pa_t v\gtrsim-C_0-\|u\|_{L^{\iy}_T(L^1_{\om})}\text{ and }
\BM^-_{\fL_0}v-\pa_t v\lesssim C_0+\|u\|_{L^{\iy}_T(L^1_{\om})}\text{
in $Q_1$. }$$ Hence we complete the proof by applying Theorem 5.2
\cite{KL} to $v$. \qed

In the following lemma, we get a useful estimate which can be
derived from Morrey's inequality.

\begin{lemma} If $\,u\in L^{\iy}_T(L^1_{\om})$ is a function with
$u(\cdot,t)\in C^1(B_r)$ for $t\in(-r^{\sm},0]$ and
$[u]_{C^{1,\ap}(Q_r)}<\iy$ for some $\ap\in(0,1)$, then we have that
$$[u]_{C^{0,1}_x(Q_r)}\le\f{(2^{\ap}+3^{\ap})\,r^{\ap}}{2^{\ap}}\,[u]_{C_x^{1,\ap}(Q_r)}
+\f{8}{\om(B_{r/2})\,r}\,\|u\|_{L^{\iy}_T(L^1_{\om})}.$$
\end{lemma}

\pf From Lemma 2.2, there is some $z_0^t\in B_r$ (depending on $t$)
such that
$$|Du(z_0^t,t)|\le\bigl(\f{3r}{2}\bigr)^{\ap}[u]_{C_x^{1,\ap}(Q_r)}
+\f{8}{\om(B_{r/2})\,r}\,\|u\|_{L^{\iy}_T(L^1_{\om})}.$$ Thus it
follows from this fact and Morrey's inequality that
\begin{equation*}\begin{split}[u]_{C_x^{0,1}(Q_r)}&\le\|u\|_{W_x^{1,\iy}(Q_r)}
\le\f{(2^{\ap}+3^{\ap})\,r^{\ap}}{2^{\ap}}\,[u]_{C_x^{1,\ap}(Q_r)}
+\f{8}{\om(B_{r/2})\,r}\,\|u\|_{L^{\iy}_T(L^1_{\om})}.
\end{split}\end{equation*}
Hence we complete the proof. \qed

\begin{thm} Let $\sm\in[\sm_0,2)$ for some $\sm_0\in(1,2)$.
Then there is some $\vr_0>0$ $($depending on $\ld,\Ld,\sm_0$ and
$n$$)$ so that if $u\in L^{\iy}_T(L^1_{\om})$ satisfies $\BI u-\pa_t u=0$ in $Q_{1+\e}$ where $\BI$ is a nonlocal, translation-invariant and
uniformly elliptic operator with respect to $\fL_*$,
then there is some $\ap>0$ such that
$$\|u\|_{C^{1,\ap}(Q_1)}\lesssim
\|u\|_{L^{\iy}_T(L^1_{\om})}.$$\end{thm}

\rk We can derive from Theorem 2.1 and Theorem 3.4 that
\begin{equation}\|u\|_{C^{\ap}(Q_1)}\lesssim\|u\|_{L^{\iy}_T(L^1_{\om})}.
\end{equation}
Also it follows from the standard telescopic sum argument \cite{CC}
and (3.5) that
\begin{equation}[u]_{C^{0,1}_x(Q_1)}\lesssim\|u\|_{L^{\iy}_T(L^1_{\om})}.\end{equation}

\pf The proof of this theorem goes along the lines of the proof of
Theorem 12.1 in \cite{CS2} by applying Theorem 3.4 to the difference
quotients in the $x$-direction
$$w^h(x,t)=\f{u(x+h,t)-u(x,t)}{|h|^{\bt}}$$ for
$\bt=\ap,2\ap,\cdots,1$. We write $w^h=w_1^h+w_2^h$ where
$w_1^h=w^h\mathbbm{1}_{Q_1}$. By Theorem 2.4 in \cite{KL}, we have that
$\,\,\BM^+_{\fL_*}w^h-\pa_t w^h\ge 0\,\,$ and $\,\,\BM^-_{\fL_*}w^h-\pa_t w^h\le
0\,\,$ in $Q_1$. Since $\pa_t w^h_2\equiv 0$ in $Q_1$ for $h$ with
$|h|<\e$, it follows from the uniform ellipticity of $\BI$ with respect to
$\fL_*$ that we have that
\begin{equation*}\begin{split} \BM^+_{\fL_0}w_1^h-\pa_t w_1^h&\ge
\BM^+_{\fL_*}w_1^h-\pa_t w_1^h\ge\BM^+_{\fL_*} w^h-\BM^+_{\fL_*}
w_2^h-\pa_t w^h\ge-\BM^+_{\fL_*}w_2^h\,\,\text{ in $Q_1$},\\
\BM^-_{\fL_0}w_1^h-\pa_t w_1^h&\le\BM^-_{\fL_*}w_1^h-\pa_t
w_1^h\le\BM^-_{\fL_*} w^h-\BM^-_{\fL_*} w_2^h-\pa_t
w^h\le-\BM^-_{\fL_*}w_2^h\,\,\text{ in
$Q_1$}.\end{split}\end{equation*} If we can show that
$|\BM^+_{\fL_*}w_2^h|\vee|\BM^-_{\fL_*}w_2^h|\lesssim
\|u\|_{L^{\iy}_T(L^1_{\om})}$ in $Q_1$, then we have that
$$\BM^+_{\fL_0}w_1^h-\pa_t w_1^h\gtrsim-c\|u\|_{L^{\iy}_T(L^1_{\om})}\text{ and
}\BM^-_{\fL_0}w_1^h-\pa_t w_1^h\lesssim
\|u\|_{L^{\iy}_T(L^1_{\om})}\text{ in $Q_1$}$$ for $h$ with
$|h|<\e$. Indeed, it can be obtained from the fact that
\begin{equation*}\begin{split}&\int_{\BR^n\s
B_{\rho}}|u(x+y,t)|\f{|K(x,y,t)-K(x,y-h,t)|}{|h|}\,dy\\
&\qquad+\int_{\BR^n\s B_{\rho}}|u(x+y+h,t)|K(x,y,t)\,dy\lesssim
\|u\|_{L^{\iy}_T(L^1_{\om})}\end{split}\end{equation*} for some
$\rho>0$ (this can be seen by using (2.8) and (3.2)). Hence, by
(2.3) and (3.5), $u$ admits the required $C_x^{1,\ap}$-estimates on
$Q_1$; more precisely,
\begin{equation}\|u\|_{C_x^{1,\ap}(Q_1)}\lesssim
\|u\|_{L^{\iy}_T(L^1_{\om})}.\end{equation}

Now we are going to show that $u$ is
$C_t^{\f{1+\alpha}{\sigma}}$-H\"older continuous in $Q_1$, following
Lemma 2 in \cite{CW}. For $(x_0,t_0)\in Q_1$, we consider
$$w(x,t)=\frac{u(rx+x_0,r^{\sigma}t+t_0)-u(x_0,t_0)-r\n
u(x_0,t_0)\cdot x}{r^{1+\alpha}}$$ for any sufficiently small $r>0$.
Then $w$ solves the given parabolic equation.

Without loss of generality, by (3.6) let us assume that
$0<[u]_{C^{0,1}_x(Q_1)}<\iy$.  Then $C_x^{1,\alpha}$-regularity of
$u$ on $Q_1$ and Lemma 3.5 imply the estimate
\begin{equation}0<[u]_{C^{0,1}_x(Q_1)}\le 5\,[u]_{C_x^{1,\ap}(Q_1)}
+\f{8}{\om(B_{1/2})}\,\|u\|_{L^{\iy}_T(L^1_{\om})}<\iy.
\end{equation}
So, by dividing $u$ by the right-hand side in the above, we assume
that $[u]_{C^{0,1}_x(Q_1)}\le 1$. We consider the function $\phi\,$
given by
\begin{equation*}\phi(y)=\begin{cases} |y|^2, &y\in B_{1+\e}, \\
                                       (1+\e)^2, & y\in\BR^n\s
                                       B_{1+\e}.
\end{cases}\end{equation*}
If $x\in B_1$, then we have that
\begin{equation*}\begin{split}L\phi(x)&=\int_{|y|<\e}[\phi(x+y)+\phi(x-y)-2\phi(x)]K(x,y,t)\,dy\\
&\qquad+\int_{|y|\ge\e}[\phi(x+y)+\phi(x-y)-2\phi(x)]K(x,y,t)\,dy\\
&\le 2\Ld\om_n\e^{2-\sm}+\f{4(2-\sm)}{\sm}\Ld\om_n\e^{-\sm}\le
6\Ld\om_n\e^{-\sm}
\end{split}\end{equation*}
for any $L\in\fL_*$, and we have that $\BI\phi\le
6\Ld\om_n\e^{-\sm}$ on $B_1$. 
Set $M=\sup_{B_1\times(-1,t_1)}w$. Then we may assume that $M\ge 0$; otherwise, we could use $-w$ in place of $w$.
If $M=w(x_M,t_M)$ for some $x_M\in B_1$ and
$t_M\in(-1,t_1)$ where $t_1=-1+\f{1}{12\Ld\om_n\e^{-\sm}}$, then it
is easy to check that the functions
\begin{equation*}\begin{split}
\phi_1(x,t)&=M\biggl(t+1+\f{[u]_{C^{0,1}_x(Q_1)}\,\phi(x)}{12\Ld\om_n\e^{-\sm}(1+\e)^2}\biggr),\\
\phi_2(x,t)&=6\Ld\om_n\e^{-\sm}M\biggl(t+1+\f{\phi(x)-\phi(x_M)}{6\Ld\om_n\e^{-\sm}}\biggr)
+[u]_{C^{0,1}_x(Q_1)},
\end{split}\end{equation*} are supersolutions of the given equation
on $Q_1$. So it follows from comparison principle \cite{KL} that
$w(x,t)\le\phi_1(x,t)\wedge\phi_2(x,t)$ for any $(x,t)\in Q_1$. We
now claim that (a) $w\le M+[u]_{C^{0,1}_x(Q_1)}$ on $Q_1$ and (b)
$M\le 4[u]_{C^{0,1}_x(Q_1)}$. Indeed, since $w\le\phi_1$ on $Q_1$
and we could assume that $6\Ld\om_n\e^{-\sm}>M\vee 1$ by the smallness of
$\e$, we can easily derive (a). For the proof of (b), if we suppose
that $M>4[u]_{C^{0,1}_x(Q_1)}$, then the inequality $w\le\phi_2$ on
$Q_1$ implies that
$$w(x_M,t_M)\le\f{1}{2}M+[u]_{C^{0,1}_x(Q_1)}<M-[u]_{C^{0,1}_x(Q_1)},$$
which is a contradiction. Hence by (3.7) and (3.8) we can get that
$$w\le 5[u]_{C^{0,1}_x(Q_1)}\le
25[u]_{C^{1,\ap}_x(Q_1)}+\f{40}{\om(B_{1/2})}\,\|u\|_{L^{\iy}_T(L^1_{\om})}
\lesssim\|u\|_{L^{\iy}_T(L^1_{\om})}$$ on $Q_1$. In a similar way, we can show that
$w\gtrsim-\|u\|_{L^{\iy}_T(L^1_{\om})}$ by constructing subsolutions
corresponding to $\phi_1$ and $\phi_2$. Thus we obtain the estimate
\begin{equation}\|u\|_{C_t^{\f{1+\ap}{\sm}}(Q_1)}\lesssim
\|u\|_{L^{\iy}_T(L^1_{\om})}\end{equation} Therefore by (3.7) and
(3.9) we obtain the required estimate. \qed

\begin{thm} Let $\bI$ be the nonlocal operator as in $(1.1)$. Then
the operator $\bI$ satisfies the following properties;

$(a)$ $\bI u(x,t)$ is well-defined for any $u\in C^{1,1}_x(x,t)\cap
L^{\iy}_T(L^1_{\om})$,

$(b)$ $\bI u$ is continuous in $\Om_{\tau}$, whenever $u\in
C^{1,1}_x(\Om_{\tau})\cap L^{\iy}_T(L^1_{\om})$.
\end{thm}

\pf (a) It can be shown as in the elliptic case.

(b) Take any $u\in L^{\iy}_T(L^1_{\om})$ and $\vep>0$. Then for any
$t\in[-\tau,0]\subset(-T,0]$, there is some $g_t\in
C_c^{\iy}(\BR^n)$ with $\supp(g_t)\supset\Om$ such that
$\|g_t-u(\cdot,t)\|_{L^1_{\om}}<\vep$. We consider a function $g\in
L_c^{\iy}(\BR^n_T)$ ( i.e. $g\in L^{\iy}(\BR^n_T)$ with compact
support ) so that
\begin{equation*}g(x,t)=\begin{cases}g_t(x), &
(x,t)\in\BR^n\times[-\tau,0],\\ 0, & (x,t)\in\BR^n\times(-T,-\tau-1)
\end{cases}\end{equation*} and $\sup_{t\in[-\tau,0]}\|g(\cdot,t)-u(\cdot,t)\|_{L^1_{\om}}<\vep$.
So we may assume that $u\in L_c^{\iy}(\BR^n_T)\cap
C_x^{1,1}(\Om_{\tau})\cap L^{\iy}_T(L^1_{\om})$. Thus it easily
follows from the continuity of $u(\cdot,t)$ in time variable on the
norm $\|\cdot\|_{L^1_{\om}}$ and the proof of the elliptic case (see
\cite{CS1}, \cite{CS2} and \cite{KL}). \qed

\section{Boundary estimates and Global estimates}

In this section, we realize that a modulus of continuity on the
parabolic boundary of the domain of some equation makes it possible
to obtain another modulus of continuity inside the domain. This can
be established by controlling the growth of $u$ away from its
parabolic boundary values via barriers, scaling and interior
regularity.

We use a barrier function which was used in \cite{CS1} for the
elliptic case and adapted to our parabolic setting. This barrier
function is appropriate as a supersolution of
$\BM^+_{\sm}\psi-\pa_t\psi\le 0$ for all values of $\sm$ greater
than a given $\sm_0$, where $\BM^+_{\sm}$ denotes the maximal
operator $\BM^+_{\fL_0(\sm)}$. Another way to say this would be to
define a larger class $\fL$ which is the union of all classes
$\fL_0(\sm)$ for $\sm\in(\sm_0,2)$, then
$\BM^+_{\fL}\psi-\pa_t\psi\le 0$. The proof of the following lemma
can be achieved by a little modification to our parabolic setting
(refer to \cite{CS1}), and so we leave the proof for the reader.

\begin{lemma} Let $\sm_0\in(0,2)$ be given. Then, for any $\sm\in(\sm_0,2)$ and $\gm\in(0,1)$,
there are some $\ap>0$ and $r>0$ so small that the function
$g_{\ap}(x,t)=(|x|-1)_+^{\ap}$ satisfies $\BM^+_{\sm}
g_{\ap}\le-1/[(2^{\sm}-1)\gm^{\sm}]$ in $(B_{1+r}\s
B_1)\times(-T,0]$.\end{lemma}

\begin{cor} Let $\sm_0\in(0,2)$ be given. Then, for any $\sm\in(\sm_0,2)$ and $\gm\in(0,1)$,
there is a continuous
function $\psi$ defined on $\BR^n_T$ such that $(a)$ $\psi=0$ in
$Q_1$, $(b)$ $\psi\ge 0$ in $\BR^n_T$, $(c)$ $\psi\ge\gm^{-\sm}$ in
$\BR^n_T\s Q_2$, $(d)$ $\BM^+_{\sm}\psi-\pa_t\psi\le 0$ and
$\pa_t\psi\ge-[(2^{\sm}-1)\gm^{\sm}]^{-1}$ in $\BR^n_T\s
Q_1$.\end{cor}

\pf We consider
$\psi(x,t)=\min\{\gm^{-\sm},C(|x|-1)_+^{\ap}\}+(t+1)_-/[(2^{\sm}-1)\gm^{\sm}]$
for some large constant $C>0$ and apply Lemma 4.1. \qed

The function $\psi$ obtained in Corollary 4.2 shall be utilized as a
barrier to prove the boundary continuity of solutions to nonlocal
parabolic equations. We observe that $\psi$ is a supersolution
outside the parabolic cube $Q_1$.

\begin{thm} Let $\sm\in(\sm_0,2)$ for $\sm_0\in(1,2)$. If
$u\in L^{\iy}_T(L^1_{\om})$ satisfies that
$$\BM^+_{\sm}u-\pa_t u\ge-C\,\,\text{ and }\,\,\BM^-_{\sm}u-\pa_t u\le C\,\,\text{
in $Q_{1+\e}$,}$$
$$|u(y,s)-u(x,t)|\le\rho((|x-y|^{\sm}+|t-s|)^{1/\sm})$$
for every $(x,t)\in\pa_p Q_1$ and $(y,s)\in\BR^n_T\s Q_1$, where
$\rho$ is a modulus of continuity, then there is another modulus of
continuity $\bar\rho$ $($depending only on
$\rho,\ld,\Ld,\sm_0,n,\|u\|_{L^{\iy}_T(L^1_{\om})}$, $\e$ and $C$,
but not on $\sm$$)$ such that
$$|u(y,s)-u(x,t)|\le\bar\rho((|x-y|^{\sm}+|t-s|)^{1/\sm})$$ for
every $(x,t)\in\overline Q_1$ and $(y,s)\in\BR^n_T$.\end{thm}

\begin{lemma} Let $\sm\in(\sm_0,2)$ for $\sm_0\in(1,2)$. If
$u\in L^{\iy}_T(L^1_{\om})$ is a function such that
$$\BM^+_{\sm}u-\pa_t u\ge-C\,\,\text{ in $Q_{1+\e}$,}$$
$$u(x,t)-u(x_0,t_0)\le\rho((|x-x_0|^{\sm}+|t-t_0|)^{1/\sm})$$
for every $(x_0,t_0)\in\pa_p Q_1$ and $(x,t)\in\BR^n_T\s Q_1$, where
$\rho$ is a modulus of continuity, then there is another modulus of
continuity $\tilde\rho$ $($depending only on
$\rho,\ld,\Ld,\sm_0,n,\|u\|_{L^{\iy}_T(L^1_{\om})},$ $\e$ and $C$,
but not on $\sm$$)$ such that
$$u(x,t)-u(x_0,t_0)\le\tilde\rho((|x-x_0|^{\sm}+|t-t_0|)^{1/\sm})$$ for
every $(x_0,t_0)\in\pa_p Q_1$ and $(x,t)\in\BR^n_T$.\end{lemma}

\pf If we write $v=u\mathbbm{1}_{Q_{1+\e}}$, then as before we have
that
\begin{equation*}\BM^+_{\sm}v-\pa_t v\ge-c_{\e}-\|u\|_{L^{\iy}_T(L^1_{\om})}\,\,\text{ and }
\,\,\BM^-_{\sm}v-\pa_t v\le
c_{\e}+\|u\|_{L^{\iy}_T(L^1_{\om})}\,\,\text{ in $Q_1$.}
\end{equation*} Since $u$ is continuous on $\overline Q_{1+\e}$, we may assume that $u\in\rB(\BR^n_T)$.

Since $\sm\ge\sm_0>0$, the function
$$p_0(x,t)=\f{1}{4}(0\vee(4-|x|^2))+(C+4\Ld\om_n(2-\sm)\f{1-3^{-\sm}}{\sm})t$$
satisfies
$\BM^+_{\sm}p_0\le-\ld\om_n+4\Ld\om_n(2-\sm)\f{1-3^{-\sm}}{\sm}$ in
$Q_1$, because
\begin{equation*}\begin{split}L p_0(x,t)&=-\int_{B_1}|y|^2 K(x,y,t)\,dy+\int_{B_3\s
B_1}\mu_t(p_0,x,y)K(x,y,t)\,dy\\
&\le-\f{\ld\om_n}{2}+4\Ld\om_n(2-\sm)\f{1-3^{-\sm}}{\sm}
\end{split}\end{equation*} for any $L\in\fL_0(\sm)$ and all $(x,t)\in Q_1$.
Since $\pa_t p_0=C+4\Ld\om_n(2-\sm)\f{1-3^{-\sm}}{\sm}$ in $Q_1$, we
have that
\begin{equation}\BM^+_{\sm}\bigl(u-p_0\bigr)-\pa_t\bigl(u-p_0\bigr)\ge
\BM^+_{\sm}u-\pa_t u-\BM^+_{\sm}p_0+\pa_t p_0\ge\f{\ld\om_n}{2}\ge
0\,\,\text{ in $Q_1$.}\end{equation} Let $\rho_0$ be the modulus of
continuity of the function $\psi$ in Corollary 4.2 and let $\rho_1$
be the modulus of continuity of the function $p$. By the assumption,
we see that
\begin{equation}\begin{split}&u(x,t)-p_0(x,t)-u(x_0,t_0)+p_0(x_0,t_0)\\
&\qquad\le\rho((|x-x_0|^{\sm}+|t-t_0|)^{1/\sm})+\rho_1((|x-x_0|^{\sm}+|t-t_0|)^{1/\sm})\end{split}\end{equation}
for every $(x_0,t_0)\in\pa_p Q_1$ and $(x,t)\in\BR^n_T\s Q_1$.

Fix any $(x_0,t_0)\in\pa_p Q_1$. For $r>0$, we define $\bar\rho$ by
$$\bar\rho(r)=\inf_{\gm\in(0,1)}\biggl(\rho(3\gm)+\rho_1(3\gm)+\bigl\|u-p_0-u(x_0,t_0)+
p_0(x_0,t_0)\bigr\|_{L^{\iy}}\rho_0\biggl(\f{r}{\gm}\biggr)\biggr).$$
Then we must show that $\bar\rho$ is a modulus of continuity. We
easily see that $\bar\rho$ is clearly monotonically increasing
because $\rho_0$ is. So we have only to show that for any $\vep>0$,
there is some $r>0$ such that $\bar\rho(r)<\vep$. Indeed, we choose
some $\gm\in(0,1)$ such that $\rho(3\gm)+\rho_1(3\gm)<\vep/2$, and
then choose some $r>0$ so that $\|u-p_0-u(x_0,t_0)+
p_0(x_0,t_0)\|_{L^{\iy}}\,\rho_0(r/\gm)<\vep/2$. Finally, we show
that there is a modulus of continuity $\tilde\rho$ such that
$u(x,t)-u(x_0,t_0)\le\tilde\rho((|x-x_0|^{\sm}+|t-t_0|)^{1/\sm})$
for any $(x,t)\in\BR^n_T$. Take any $(x,t)\in\BR^n_T$. For $\gm>0$,
we consider a barrier function
\begin{equation*}\begin{split}B(x,t)&=u(x_0,t_0)-p_0(x_0,t_0)+\rho(3\gm)+\rho_1(3\gm)\\
&\quad+\gm^{\sm}\bigl\|u-p_0-u(x_0,t_0)+
p_0(x_0,t_0)\bigr\|_{L^{\iy}}\psi\biggl(-x_0+\f{x-x_0}{\gm},\f{t-t_0}{\gm^{\sm}}\biggr).
\end{split}\end{equation*}
By (4.2) and the definition of $\psi$, we have that
$$B(x,t)\ge u(x_0,t_0)-p_0(x_0,t_0)+\rho(3\gm)+\rho_1(3\gm)\ge
u(x,t)-p_0(x,t)$$ for any $(x,t)\in Q_{3\gm}(x_0,t_0)\cap(\BR^n_T\s
Q_1)$. Also by the definition of $\psi$, we obtain that $B(x,t)\ge
u(x,t)-p_0(x,t)$ for any $(x,t)\in (\BR^n_T\s
Q_{3\gm}(x_0,t_0))\cap(\BR^n_T\s Q_1)$. Thus we have that $B\ge
u-p_0$ on $\BR^n_T\s Q_1$. By (d) of Corollary 4.2, we see that
$\BM^+_{\sm}\psi-\pa_t\psi\le 0$ in $Q_{1/\gm}(\f{\gm+1}{\gm}x_0,0)$
because $Q_{1/\gm}(\f{\gm+1}{\gm}x_0,0)\subset\BR^n_T\s Q_1$. We
observe that
\begin{equation*}\BM^+_{\sm}\psi-\pa_t\psi\le 0\,\text{ in
$Q_{1/\gm}\biggl(\f{\gm+1}{\gm}x_0,0\biggr)$}\,\,\Leftrightarrow\,\,
\BM^+_{\sm}B-\pa_t B\le 0\,\text{ in $Q_1,$}\end{equation*} by
Corollary 4.2. Taking the infimum on $\gm$, it follows from
comparison principle (Theorem 2.3 in \cite{KL}) that
$$u(x,t)-p_0(x,t)\le B(x,t)\le u(x_0,t_0)-p_0(x_0,t_0)+\bar\rho((|x-x_0|^{\sm}+|t-t_0|)^{1/\sm})$$
for all $(x,t)\in\BR^n_T$, because
$$\psi(x,t)\le|\psi(x,t)-\psi(x_0,t_0)|\le\rho_0((|x-x_0|^{\sm}+|t-t_0|)^{1/\sm}),\forall
(x,t)\in\BR^n_T.$$ Hence we have that
$u(x,t)-u(x_0,t_0)\le\tilde\rho((|x-x_0|^{\sm}+|t-t_0|)^{1/\sm})$
for all $(x,t)\in\BR^n_T$, where $\tilde\rho=\rho_1+\bar\rho$.
Therefore we complete the proof.\qed

\begin{lemma} Let $\sm\in(\sm_0,2)$ for $\sm_0\in(1,2)$. If
$u\in L^{\iy}_T(L^1_{\om})$ satisfies that
$$\BM^+_{\sm}u-\pa_t u\ge-C\,\,\text{ and }\,\,\BM^-_{\sm}u-\pa_t u\le C\,\,\text{
in $Q_{1+\e}$,}$$
$$|u(y,s)-u(x,t)|\le\rho((|x-y|^{\sm}+|t-s|)^{1/\sm})$$
for every $(x,t)\in\pa_p Q_1$ and $(y,s)\in\BR^n_T$, where $\rho$ is
a modulus of continuity, then there is another modulus of continuity
$\bar\rho$ $($depending only on
$\rho,\ld,\Ld,\sm_0,n,\|u\|_{L^{\iy}_T(L^1_{\om})},$ $\e$ and $C$,
but not on $\sm$$)$ such that
$$|u(y,s)-u(x,t)|\le\bar\rho((|x-y|^{\sm}+|t-s|)^{1/\sm})$$ for
every $(x,t)\in\overline Q_1$ and $(y,s)\in\BR^n_T$.\end{lemma}

\pf If we set $v=u\mathbbm{1}_{Q_{1+\e}}$, then as before we have
that
\begin{equation*}\BM^+_{\sm}v-\pa_t v\ge-c_{\e}-\|u\|_{L^{\iy}_T(L^1_{\om})}\,\,\text{ and }
\,\,\BM^-_{\sm}v-\pa_t v\le
c_{\e}+\|u\|_{L^{\iy}_T(L^1_{\om})}\,\,\text{ in $Q_1$.}
\end{equation*} Since $u$ is continuous on $\overline Q_{1+\e}$, we may
assume that $u\in\rB(\BR^n_T)$. Hence it can easily be obtained by
an adaptation of Lemma 3 in [CS1] to our parabolic setting. \qed

{\bf Proof of Theorem 4.3.} We apply Lemma 4.4 to both $u$ and $-u$
to obtain a modulus of continuity that applies from any point on
$\pa_p Q_1$ to any point in $\BR^n_T$. Then we use Lemma 4.5 to
finish the proof. \qed

\section{Some results by approximation}

In this section, we show that two equations which are very close to
each other in some appropriate way have their solutions which are
close by each other on the unit cube $Q_1$.

In what follows, for a function $u:\BR^n\to\BR$ and a parabolic
quadratic polynomial $p$ we denote by $u^p_{Q_r(x_0,t_0)}\fd
p\mathbbm{1}_{Q_r(x_0,t_0)}+u\mathbbm{1}_{\BR^n_T\s Q_r(x_0,t_0)}$.
The following lemma is an usual result in analysis on viscosity
solutions, and so we will skip the proof.

\begin{lemma} Let $\BI$ be a uniformly elliptic operator in the
sense of $(1.2)$ with respect to some class $\fL$ and let
$u:\BR^n_T\to\BR$ be a function which is upper semicontinuous on
$\overline\Om_{\tau}$. Then the followings are equivalent.

$(a)$ $u$ is a viscosity subsolution of $\,\BI u-\pa_t u=f$ in
$\Om_{\tau}$, i.e. $\BI u-\pa_t u\ge f$ in $\Om_{\tau}$.

$(b)$ If $p$ is a parabolic quadratic polynomial satisfying
$u(x_0,t_0)=p(x_0,t_0)$ and $u\le p$ in $Q_r(x_0,t_0)$ where
$Q_r(x_0,t_0)\subset\Om_{\tau}$ for some $r>0$ and
$(x_0,t_0)\in\Om_{\tau}$, then we have that $\BI
u^p_r(x_0,t_0)-\pa_t u^p_r(x_0,t_0)\ge f(x_0,t_0)$ for $u^p_r=
u^p_{Q_r(x_0,t_0)}$.\end{lemma}

We want to show that if $\BI_k u_k(x,t)=f_k(x,t)$ and $\BI_k\to\BI$,
$u_k\to u$ and $f_k\to f$ in some appropriate way, then $\BI
u(x,t)=f(x,t)$.

In the elliptic case \cite{CS1}, the solution space $L^1_{\om}$ is
enough for the weakly convergence of operators $\BI_k$. In the
parabolic case, the possible substitute for the solution space
$L^1_{\om}$ is $L^{\iy}_T(L^1_{\om})$. This makes it possible to
obtain the stability properties for the nonlocal parabolic case.

\begin{definition} We say that $\BI_k$ {\rm converges weakly to $\BI$} in $\Om_{\tau}$
with respect to $\om$ $($and we denote by
$\lim_{k\to\iy}\BI_k=_{\om}\BI$ in $\Om_{\tau}$$)$, if for any
$(x_0,t_0)\in\Om_{\tau}$ there is some
$Q_r(x_0,t_0)\subset\Om_{\tau}$ such that
\begin{equation}\lim_{k\to\iy}\BI_k u^p_r=\BI u^p_r\end{equation} uniformly in
$Q_{r/2}(x_0,t_0)$  for any function $u_r^p$ of the form
$u^p_r=u^p_{Q_r(x_0,t_0)}$ where $p$ is a parabolic quadratic
polynomial and $u\in L^{\iy}_T(L^1_{\om})$.
\end{definition}

\begin{lemma} Let $\BI$ be a uniformly elliptic operator with
respect to a class $\fL$ of linear integro-differential operators.
If $u^p_r=u^p_{Q_r(x_0,t_0)}$ where $p$ is a parabolic quadratic
polynomial and $u\in L^{\iy}_T(L^1_{\om})$, then $\BI u^p_r$ is
continuous in $Q_r(x_0,t_0)$.
\end{lemma}

\pf Since $u^p_r\in C_x^{1,1}(Q_r(x_0,t_0))\cap
L^{\iy}_T(L^1_{\om})$, the required result easily follows from
Theorem 3.7. \qed

\begin{lemma} Let $\{\BI_k\}$ be a sequence of uniformly elliptic
operators with respect to some class $\fL$. Assume that Assumption
3.3 holds. Let $\{u_k\}\subset L^{\iy}_{\tau}(L^1_{\om})$ be a
sequence of lower semicontinuous functions in $\Om_{\tau}$ such that

$(a)$ $\BI_k u_k-\pa_t u_k\le f_k$ in $\Om_{\tau}$, $(b)$
$\lim_{k\to\iy}u_k=u$ in the $\Gm$ sense in $\Om_{\tau}$,

$(c)$ $\lim_{k\to\iy}\|u_k-u\|_{L^{\iy}_{\tau}(L^1_{\om})}=0$, $(d)$
$\lim_{k\to\iy}\BI_k=_{\om}\BI$ in $\Om_{\tau}$,

$(e)$ $\lim_{k\to\iy}f_k=f$ locally uniformly in $\Om_{\tau}$,

$(f)$ $\sup_{k\in\BN}\sup_{\Om_{\tau}}|u_k|\le C<\iy$.

\noindent Then we have that $\BI u-\pa_t u\le f$ in
$\Om_{\tau}$.\end{lemma}

\pf Let $p$ be a parabolic quadratic polynomial touching $u$ from
below at a point $(x,t)$ in a neighborhood $V\subset\Om_{\tau}$.
Since $\{u_k\}$ $\Gm$-converges to $u$ in $\Om_{\tau}$, there are a
cube $Q_r(x,t)\subset V$ and a sequence $\{(x_k,t_k)\}\subset
Q_r(x,t)$ with $\lim_{k\to\iy}\dd((x_k,t_k),(x,t))=0$ such that $p$
touches $u_k$ from below at $(x_k,t_k)$(refer to \cite{GD}). Without
loss of generality, we assume that $Q_r(x,t)$ is a cube so that
(5.1) holds for the point $(x,t)$.

If $(u_k)^p_r=(u_k)^p_{Q_r(x,t)}$, then we have $\BI
(u_k)^p_r(x_k,t_k)-\pa_t(u_k)^p_r(x_k,t_k)\le f(x_k,t_k)$. If we set
$u^p_r=u^p_{Q_r(x,t)}$, then we see that $u^p_r(z,q)=(u_k)^p_r(z,q)$
and $\pa_t(u_k)^p_r(z,q)=\pa_t u^p_r(z,q)$ for any $k\in\BN$ and
$(z,q)\in Q_{r}(x,t)$. Take any $(z,q)\in Q_{r/4}(x,t)$. Then we
have that
\begin{equation*}\begin{split}&|\BI_k (u_k)^p_r(z,q)-\pa_t(u_k)^p_r(z,q)-\BI u^p_r(z,q)+\pa_t
u^p_r(z,q)|\\
&\le|\BI_k (u_k)^p_r(z,q)-\BI_k u^p_r(z,q)|+|\BI_k
u^p_r(z,q)-\BI u^p_r(z,q)|\\
&\le|\BM^+_{\fL}((u_k)^p_r-u^p_r)(z,q)|\vee|\BM^+_{\fL}(u^p_r-(u_k)^p_r)(z,q)|+|\BI_k
u^p_r(z,q)-\BI u^p_r(z,q)|\\
&\le\sup_{L\in\fL}|L((u_k)^p_r-u^p_r)(z,q)|+|\BI_k
u^p_r(z,q)-\BI u^p_r(z,q)|\\
&\le\int_{\BR^n\s
B_{r/2}}|\mu_q((u_k)^p_r-u^p_r,z,y)|K(x,y,t)\,dy+|\BI_k
u^p_r(z,q)-\BI u^p_r(z,q)|\\ &\le C_r\int_{\BR^n\s
B_{r/2}}\bigl\{|((u_k)^p_r-u^p_r)(z+y,q)|+|((u_k)^p_r-u^p_r)(z-y,q)|\bigr\}\om(y)\,dy\\
&\qquad\qquad\qquad\qquad\qquad\qquad\qquad\qquad\qquad+|\BI_k
u^p_r(z,q)-\BI u^p_r(z,q)|\\
&\le C\int_{\BR^n}2|(u_k)^p_r(y,q)-u^p_r(y,q)|\sup_{z\in
B_{r/4}}\om(y+z)\,dy+|\BI_k
u^p_r(z,q)-\BI u^p_r(z,q)|\\
&\le C\|u_k-u\|_{L^{\iy}_{\tau}(L^1_{\om})} +|\BI_k u^p_r(z,q)-\BI
u^p_r(z,q)|\to 0\end{split}\end{equation*} as $k\to\iy$, by using
Assumption 3.3 and (3.2). Since $(u_k)^p_r\in C^2(Q_r(x,t))\cap
L^{\iy}_{\tau}(L^1_{\om})$ for all $k\in\BN$ and
$\lim_{k\to\iy}\|u_k-u\|_{L^{\iy}_{\tau}(L^1_{\om})}=0$, we see
$u^p_r\in C^2(Q_r(x,t))\cap L^{\iy}_{\tau}(L^1_{\om})$, and thus
$\BI u^p_r$ is continuous in $Q_r(x,t)$ (by Lemma 5.3). Thus by
(5.1) we have that
\begin{equation*}\begin{split}&|\BI_k (u_k)^p_r(x_k,t_k)-\pa_t(u_k)^p_r(x_k,t_k)-\BI u^p_r(x,t)+\pa_t
u^p_r(x,t)|\\&\qquad\qquad\le|\BI_k
(u_k)^p_r(x_k,t_k)-\pa_t(u_k)^p_r(x_k,t_k)-\BI u^p_r(x_k,t_k)+\pa_t
u^p_r(x_k,t_k)|\\
&\qquad\qquad\quad+|\BI_k u^p_r(x_k,t_k)-\BI u^p_r(x_k,t_k)|+|\BI u^p_r(x_k,t_k)-\BI u^p_r(x,t)|\\
&\qquad\qquad\quad+|\pa_t u^p_r(x_k,t_k)-\pa_t u^p_r(x,t)|\to
0\,\,\,\,\,\text{ as $k\to\iy$. }\end{split}\end{equation*} Since
$\lim_{k\to\iy}\dd((x_k,t_k),(x,t))=0$ and $\lim_{k\to\iy}f_k=f$
locally uniformly in $\Om_{\tau}$, we have that $f_k(x_k,t_k)\to
f(x,t)$. Thus this implies that $\BI u^p_r(x,t)-\pa_t u^p_r(x,t)\le
f(x,t)$. Hence we conclude that $\BI u-\pa_t u\le f$ in
$\Om_{\tau}$.\qed

\begin{lemma} Let $u^p_r=u^p_{Q_r}$ where $p$
is a quadratic polynomial and $u\in L^{\iy}_T(L^1_{\om})$. If
$\{\BI_k\}$ is a sequence of uniformly elliptic operators with
respect to some class $\fL$ satisfying Assumptions 2.2 and 2.3, then
there is a subsequence $\{\BI_{k_j}\}$ such that $\BI_{k_j} u^p_r$
converges uniformly in $Q_{r/2}$.\end{lemma}

\pf We have only to find a uniform modulus of continuity for
$\BI_{k_j} u^p_r$ in $Q_r$ so that the lemma follows from
Arzela-Ascoli Theorem.

Take any two $(x,t),(y,s)\in Q_{r/2}$ with $\dd((x,t),(y,s))<r/8$.
By the uniform ellipticity of $\BI_k$, we have that
\begin{equation}\BI_k u^p_r(x,t)-\BI_k
u^p_r(y,t)\le\BM^+_{\fL}(u^p_r-\btau_{y-x}u^p_r)(x,t).\end{equation}
Also we see that $u^p_r-\btau_{y-x}u^p_r=0$ in
$B_{r/4}(x)\times\{t\}$, because
\begin{equation*}p(x+z,t)+p(x-z,t)-2p(x,t)-p(y+z,t)-p(y-z,t)+2p(y,t)=0\end{equation*} for
any $z\in B_{r/4}$. From (3.2) and Assumption 3.3, for any $L\in\fL$
we have that
\begin{equation*}\begin{split}
&L(u^p_r-\btau_{y-x}u^p_r)(x,t)=\int_{\BR^n\s B_{r/4}}\mu_t(u^p_r-\btau_{y-x}u^p_r,x,z)K(x,z,t)\,dz\\
&=\int_{\BR^n\s B_{r/4}}\bigl[u^p_r(x+z,t)+u^p_r(x-z,t)-2
u^p_r(x,t)\bigr]\,K(x,z,t)\,dz\\
&\qquad-\int_{\BR^n\s B_{r/4}}\bigl[u^p_r(y+z,t)+u^p_r(y-z,t)-2 u^p_r(y,t)\bigr]\,K(x,z,t)\,dz\\
&\le C_r\bigl(|p(x,t)-p(y,t)|+\|\btau_x u^p_r-\btau_y
u^p_r\|_{L^{\iy}_T(L^1_{\om})}\bigr)\\
&\le C_r\sup_{\substack{|\xi-\e|\le|x-y|\\\xi,\e\in
B_{r/2}}}\biggl(|p(\xi,t)-p(\e,t)|+\int_{\BR^n}|\btau_{\xi-\e}u^p_r(z,t)-u^p_r(z,t)|\bigl(\sup_{B_r(z)}\om\bigr)\,dz\biggr)\\
&\le\fm_1(|x-y|)
\end{split}\end{equation*} where $\fm_1$ is defined as
$\fm_1(\vr)=\sup_{t\in(-r/2,r/2]}\fm_t(\vr)$ and
$$\fm_t(\vr)=C_r\sup_{\substack{|\xi-\e|\le\vr\\\xi,\e\in
B_{r/2}}}\biggl(|p(\xi,t)-p(\e,t)|+\int_{\BR^n}|\btau_{\xi-\e}u^p_r(\,\cdot\,,t)-u^p_r(\,\cdot\,,t)|\,\om(z)\,dz\biggr).$$
Thus by (5.2) we obtain that
\begin{equation}\begin{split}&\BI_k u^p_r(x,t)-\BI_k
u^p_r(y,t)\le\fm_1(|x-y|).\end{split}\end{equation}

On the other hand, we now estimate $\BI_k u^p_r(y,t)-\BI_k
u^p_r(y,s)$. By the uniform ellipticity of $\bI_k$, we have that
\begin{equation*}\BI_k u^p_r(y,t)-\BI_k u^p_r(y,s)\le\bM^+_{\fL}(u_r^p-\btau^{s-t}u_r^p)(x,t).\end{equation*}
Observing $u^p_r-\btau^{s-t}u^p_r=0$ in $B_{r/4}(x)\times\{t\}$ from
the fact that
\begin{equation*}p(x+z,t)+p(x-z,t)-2p(x,t)-p(x+z,s)-p(x-z,s)+2p(x,s)=0\end{equation*}
for any $z\in B_{r/4}$, as in the above we can obtain that
\begin{equation}\BI_k u^p_r(y,t)-\BI_k u^p_r(y,s)\le\fm_2(|t-s|^{1/\sm})\end{equation}
for some modulus of continuity $\fm_2$ depending on $u$ but not on
$\BI_k$. Hence by (5.3) and (5.4) we conclude that
\begin{equation*}\begin{split} &\BI_k u^p_r(x,t)-\BI_k
u^p_r(y,s)\le\fm((|x-y|^{\sm}+|t-s|)^{1/\sm})
\end{split}\end{equation*} where
$\fm((|x-y|^{\sm}+|t-s|)^{1/\sm})=\fm_1(|x-y|)+\fm_2(|t-s|^{1/\sm})$.
Here it is clear that $\fm(\vr)$ is a modulus of continuity
depending on $u$ but not on $\BI_k$. Therefore there is a
subsequence that converges uniformly by Arzela-Ascoli Theorem. \qed

\begin{thm} Let $\{\BI_k\}$ be a sequence of uniformly elliptic
operators with respect to some class $\fL$ satisfying Assumptions
3.2 and 3.3 Then there is a subsequence $\{\BI_{k_j}\}$ that
converges weakly. \end{thm}

\pf Since the space $L^{\iy}_T(L^1_{\om})$ is separable with respect
to the norm $\|\cdot\|_{L^{\iy}_T(L^1_{\om})}$, we can take a
countable dense subset $\cD:=\{u_j\}$ of $L^{\iy}_T(L^1_{\om})$. We
note that the set $\Pi_{n+1}$ of all parabolic quadratic polynomials
is a finite dimensional space which has a countable dense subset
$\{p_j\}$. For each $k\in\BN$, we set
$$v_{k,j_1,j_2}=p_{j_1}\mathbbm{1}_{Q_{2^{-k}}}+u_{j_2}\mathbbm{1}_{\BR^n_T\s
Q_{2^{-k}}}.$$ Take any $\vep>0$ and any $v$ as in (5.1), i.e.
$v=u^p_{Q_r}$ for some $p\in\Pi_{n+1}$ and $r>0$. Then we choose $k$
so that $2^{-k}<r<2^{-k+1}$ and select some $j_1$ and $j_2$ such
that $\|u_{j_2}-u\|_{L^{\iy}_T(L^1_{\om})}<\vep$, and $|D_x^2
p_{j_1}-D_x^2 p|<\vep$, $|D_x p_{j_1}-D_x p|<\vep$ and
$|p_{j_1}-p|<\vep$ in $Q_{2^{-k}}$. Since the set
$\cE=\{v_{k,j_1,j_2}\}$ is countable and dense, we can arrange it in
a sequence $v_j$ of the form
$v_j=p_j\mathbbm{1}_{Q_{r_j}}+u_j\mathbbm{1}_{\BR^n_T\s Q_{r_j}}$ so
that for each $v$ as in (5.1) there is some $v_j$ such that
\begin{equation}\begin{split}\|v-v_j\|_{L^{\iy}_T(L^1_{\om})}&<\vep,\\
\sup_{Q_{r/2}}|D^k v-D^k v_j|&<\vep\,\text{for all $k=0,1,2$.
}\end{split}\end{equation} By Lemma 5.5, for each $v_j\in\cE$ there
exists a subsequence $\BI_{k_i}$ such that $\BI_{k_i}(v_j)$
converges uniformly in $Q_{r_j/2}$.

By a standard diagonalization process, there is a subsequence
$\{\BI_{k_i}\}$ such that for each $v_j$, $\{\BI_{k_i} v_j\}$
converges uniformly in $Q_{r_j/2}$. We call this limit
$\BI_{\ast}v_j (x,t)$. If $v$ is any test function, then there is
some $j$ such that $v$ is close enough to $v_j$ in the sense of
(5.2). Take any $(x,t)\in Q_{r/2}$. By the mean value theorem, we
see that
$$\mu_t(v-v_j,x,y)=\int_0^1\int_0^1\la D_x^2(v-v_j)((x+\tau
y)-2s\tau y,t) y,y\ra\,ds\,d\tau$$ for any $y\in B_{r/2}$. Thus it
follows from (3.2), (5.2), (5.5) and Assumption 3.3 that
\begin{equation*}\begin{split}\BI_k v(x,t)-\BI_k v_j(x,t)&\le\left(\int_{B_{r/2}}+\int_{\BR^n\s
B_{r/2}}\right)\mu_t(v-v_j,x,y)K_{\fL}(x,y,t)\,dy\\
&\lesssim\vep+\|v-v_j\|_{L^{\iy}_T(L^1_{\om})}\lesssim\vep.
\end{split}\end{equation*} for any $(x,t)\in Q_{r/2}$, uniformly in $k$.
Taking $i$ large enough, we thus have that
$$|\BI_{k_i}v(x,t)-\BI_{\ast}v_j(x,t)|<2\vep,$$ and thus
$\{\BI_{k_i}v(x,t)\}$ is a Cauchy sequence in $L^{\iy}(Q_{r/2})$. We
define $\BI_{\ast}v(x,t)$ to be the uniform limit of this sequence
in $Q_{r/2}$. Thus we have shown that $\{\BI_k v(x,t)\}$ converges
uniformly to $\BI_{\ast}v(x)$ in $Q_{r/2}$.

To finish the proof, we must show that the operator $\BI_{\ast}$ can
be extended to a uniformly elliptic operator for all test functions
$\vp$. We note that for any two test functions $v_1,v_2\in\cT$, we
have that
$$\BM^-_{\fL}(v_1-v_2)(x,t)\le\BI_k(v_1-v_2)(x,t)\le
\BM^+_{\fL}(v_1-v_2)(x,t).$$ Passing to the limit in this
inequality, we obtain that $$\BM^-_{\fL}(v_1-v_2)(x,t)\le
\BI_{\ast}(v_1-v_2)(x,t)\le\BM^+_{\fL}(v_1-v_2)(x,t).$$
Approximating on an arbitrary test function $\vp$ as in the proof of
Lemma 5.1, we can extend $\BI_{\ast}$ in a unique way to all test
functions $\vp$ such that $\BI_{\ast}$ is uniformly elliptic with
respect to $\fL$. \qed

\begin{lemma} For some $\sm\ge\sm_0>1$ and $\gm\in(0,\sm_0-1)$, let $\BI_0,\BI_1$ and
$\BI_2$ be nonlocal uniformly elliptic operators with respect to
$\fL_0(\sm)$ satisfying Assumptions 3.2 and 3.3. Suppose that the
boundary value problem
\begin{equation*}\begin{cases} \BI_0 u-\pa_t u=0 &\text{ in $Q_1$, }\\
u=h &\text{ in $\BR^n_T\s Q_1$ }\end{cases}\end{equation*} has at
most one solution $u$ for any $h\in L^{\iy}_T(L^1_{\om})$. Given a
modulus of continuity $\vr$ and $\vep>0$, there are a small $\dt>0$
and a large $R>0$ so that if $u,v,\BI_0,\BI_1$ and $\BI_2$ satisfy
\begin{equation*}\begin{split} \BI_0
v-\pa_t v=0,\,\BI_1 u-\pa_t u&\ge-\dt,\,\BI_2 u-\pa_t
u\le\dt,\,\|\BI_1-\BI_0\|\vee
\|\BI_2-\BI_0\|\le\dt \text{ in $Q_1$, }\\
u&=v\,\text{ in $\BR^n_T\s Q_1$, }\\
|u(x,t)-u(y,s)|&\vee|v(x,t)-v(y,s)|\le\vr((|x-y|^{\sm}+|t-s|)^{1/\sm})\\&\qquad\quad\text{
for any $x\in
Q_R\s Q_1$ and $y\in \BR^n_T\s Q_1$,}\\
|u(x,t)|&\le
M(1\vee(|x|^{\sm}+|t|)^{\f{1+\gm}{\sm}}),\end{split}\end{equation*}
then we have that $|u-v|<\vep$ in $Q_1$.\end{lemma}

\pf Assume that the result was not true. Then there would be
sequences $\{R_k\},\{\BI_0^{(k)}\},$
$\{\BI_1^{(k)}\},\{\BI_2^{(k)}\},$
$\{\dt_k\},\{u_k\},\{v_k\},\{\BI_k\}$ and $\{f_k\}$ such that
$R_k\to\iy$, $\dt_k\to 0$ and all the assumptions of the lemma hold,
but $\sup_{Q_1}|u_k-v_k|\ge\vep$. Since $\{\BI_0^{(k)}\}$ is a
sequence of uniformly elliptic operators, it follows from Theorem
5.6 that there is a subsequence that converges weakly to some
nonlocal operator $\BI_0$ which is uniformly elliptic with respect
to the same class $\fL_0(\sm)$. Moreover we see that
$\{\BI_1^{(k)}\}$ and $\{\BI_2^{(k)}\}$ converge to $\BI_0$ weakly,
because $\|\BI_0^{(k)}-\BI_1^{(k)}\|\to 0$ and
$\|\BI_0^{(k)}-\BI_2^{(k)}\|\to 0$.

Since $\vr$ is a modulus of continuity on $\pa_p Q_1$ of both
$\{u_k\}$ and $\{v_k\}$, by Theorem 4.3 there is a modulus of
continuity $\tilde\vr$ which extends to the full unit cube
$\overline Q_1$. Thus $\{u_k\}$ and $\{v_k\}$ have a modulus of
continuity on $Q_{R_k}$ with $R_k\to\iy$. We can find subsequences
$\{u_{k_j}\}$ and $\{v_{k_j}\}$ which converges uniformly on compact
sets in $\BR^n_T$ to $u$ and $v$, respectively. Since
$|u_{k_j}|\vee|v_{k_j}|\le g\in L^{\iy}_T(L^1_{\om})$ for all $j$
where $$g(x,t)=M(1\vee(|x|^{\sm}+|t|)^{\f{1+\gm}{\sm}}),$$ it
follows from the Lebesgue's dominated convergence theorem that
$u,v\in L^{\iy}_T(L^1_{\om})$ and moreover
$$\lim_{j\to\iy}\|u_{k_j}-u\|_{L^{\iy}_T(L^1_{\om})}=0\,\,\text{ and }\,\,
\lim_{j\to\iy}\|v_{k_j}-v\|_{L^{\iy}_T(L^1_{\om})}=0.$$ Since
$\sup_{Q_1}|u_{k_j}-v_{k_j}|\ge\vep$, $u$ and $v$ must be different.
By Lemma 5.4, we see that $u$ and $v$ solve the same equation $\BI_0
u-\pa_t u=\BI_0 v-\pa_t v=0$ in $Q_1$. Thus by the assumption, we
have $u=v$, which is a contradiction. \qed

\rk We will apply Lemma 5.7 to a translation-invariant operator
$\BI_0$. In case that $\BI_0$ is a translation-invariant elliptic
operator, the uniqueness for the viscosity solution of the boundary
value problem was discussed in \cite{CS2}.

We also obtain the following simplified one of Lemma 5.7. The
difference between this and Lemma 5.7 is that in Lemma 5.8 below we
fix the boundary value $h$, but we do not need a modulus of
continuity in $Q_R\s Q_1$ and also on $\pa_p Q_1$.

\begin{lemma} For some $\sm\ge\sm_0>1$,
let $\BI_0,\BI_1$ and $\BI_2$ be nonlocal uniformly elliptic
operators with respect to $\fL_0(\sm)$ satisfying Assumptions 3.2
and 3.3. Suppose that the boundary value problem
\begin{equation*}\begin{cases}\BI_0 u-\pa_t u =0 &\text{ in $Q_1$, }\\
u=h &\text{ in $\BR^n_T\s Q_1$ }\end{cases}\end{equation*} has at
most one solution $u$ for any $h\in L^{\iy}_T(L^1_{\om})$. Assume
that $h$ is continuous on $\pa_p Q_1$. Given any $\vep>0$, there is
some small $\dt>0$ so that if $u,v,\BI_0,\BI_1$ and $\BI_2$ satisfy
\begin{equation*}\begin{split} \BI_0 v-\pa_t v&=0,\,\BI_1
u-\pa_t u\ge-\dt,\,\BI_2
u-\pa_t u\le\dt,\,\|\BI_1-\BI_0\|\vee\|\BI_2-\BI_0\|\le\dt\,\text{ in $Q_1$, }\\
u&=v=h\,\,\text{ in $\BR^n_T\s Q_1$, }\end{split}\end{equation*}
then we have that $|u-v|<\vep$ in $Q_1$.\end{lemma}

\pf We proceed the proof along the same line as that of Lemma 5.7.
Assuming that the result was not true, we finish up the proof it by
getting a contradiction. Assume that there are sequences
$\{\BI_0^{(k)}\},\{\BI_1^{(k)}\},\{\BI_2^{(k)}\},$
$\{\dt_k\},\{u_k\},\{v_k\}$ and $\{\BI_k\}$ such that $\dt_k\to 0$
and all the assumptions of the lemma hold, but
$\sup_{Q_1}|u_k-v_k|\ge\vep$. The functions $u_k$ and $v_k$ have a
fixed value $h$ outside $Q_1$. Since $h$ is continuous on $\pa_p
Q_1$, by Theorem 4.3 we see that $\{u_k\}$ and $\{v_k\}$ are
equicontinuous in $\overline Q_1$. So by Arzela-Ascoli Theorem,
there is a subsequence which converges uniformly in $\overline Q_1$.

Continuing as in the proof of Lemma 5.7, we can take a subsequence
such that $\{\BI_1^{(k)}\}$ and $\{\BI_0^{(k)}\}$ converges to
$\BI_0$ weakly. Let $u$ and $v$ be the uniform limits of $u_k$ and
$v_k$ in $Q_1$, respectively. Then we have that
$\sup_{Q_1}|u_k-v_k|\ge\vep$. But by Lemma 5.4, $u$ and $v$ must
solve the same equation $\BI_0 u-\pa_t u=\BI_0 v-\pa_t v=0$ in
$Q_1$. Thus we conclude that $u=v$, which is a contradiction. \qed

\section{$C^{1,\ap}$-regularity for nonlocal parabolic equations with variable
coefficients}

The main concern of this section is to obtain $C^{1,\ap}$ estimates
for nonlocal parabolic equations which are not necessarily
translation-invariant. Since our proofs rely on rescaling argument
repeatedly, a kind of scale invariance will be needed. Even if we do
not require a particular equation to be scale invariant, we will
consider our equations within a whole class of equations that is
scale invariant for which our regularity result up to the boundary
is supposed to apply. Our proof on the parabolic case that will be
given in this section is based on that \cite{CS1} of the elliptic
case and the results \cite{KL} of the parabolic case, but the main
difference between them is to extend the solution space $\rB(\BR^n)$
on the elliptic one to the more flexible space
$L^{\iy}_T(L^1_{\om})$ on the parabolic one and to use the more
wider class of kernels involving variables $(x,t)\in\BR^n_T$.

The class $\fL$ is said to have {\it scale $\sm$} if whenever the
integro-differential operator with kernel $K(x,y,t)$ is in $\fL$,
its rescaled kernel $K_{\ld}(x,y,t):=\ld^{n+\sm}K(x,\ld y,t)$ is
also in $\fL$ for any $\ld\in(0,1)$. For example, the class $\fL_0$
defined in (1.2) has scale $\sm$, but the class $\fL_*$ defined in
(2.8) does not.

It is easy to check that if $\fL$ has scale $\sm$ and $u$ solves an
equation $\BI u(x,t)-\pa_t u(x,t)=f(x,t)$ in $Q_3$ that is elliptic
with respect to $\fL$, then the function $w_{\mu}(x,t)=\mu\,u(\ld
x,\ld^{\sm}t)$ solves a uniformly elliptic equation
\begin{equation}\BI_{\mu,\ld}w_{\mu}(x,t)-\pa_t w_{\mu}(x,t)=\ld^{\sm}\mu f(\ld
x,\ld^{\sm}t)\,\text{ in $Q_3$ }\end{equation} with respect to the
same class $\fL$. Equivalently, this condition becomes
\begin{equation}\BI_{1,\ld}w_1(x,t)-\pa_t w_1(x,t)=\ld^{\sm} f(\ld
x,\ld^{\sm}t)\,\text{ in $Q_3$; }\end{equation} that is,
$\BI_{1,\ld}u-\pa_t u=\ld^{\sm} f$ in $Q_{3\ld}$. For instance, if
$\BI u(x,t)=\int_{\BR^n}\mu_t(x,y,t)K(x,y,t)\,dy$, then
$\BI_{\mu,\ld}$ is given by
\begin{equation*}\BI_{\mu,\ld}u(x,t)=\int_{\BR^n}[u(x+y,t)+u(x-y,t)-2
u(x,t)]\ld^{n+\sm}K(x,\ld y,t)\,dy.\end{equation*} Here note that
{\it the coefficient $\mu$ does not have any effect on a linear
operator.}

Since $\fL_1$ is the largest scale invariant class contained in
the class $\fL_*$ satisfying (2.8), we observe that it follows from Theorem 3.6
that an equation $\BI u-\pa_t u=0$ has interior
$C^{1,\ap}$-estimates for some $\ap\in(0,1)$, provided that $\BI$ is
uniformly elliptic with respect to the class $\fL_1$.

Our main result in this section is to obtain that if an equation
$\BI^{(0)}u(x)-\pa_t u=0$ is uniformly elliptic with respect to a
scale invariant class with interior $C^{1,\bt}$-estimates and we
have another equation $\BI u-\pa_t u=f$ for a little perturbation
$\BI$ of $\BI^{(0)}$, then this equation also has interior
$C^{1,\ap}$-estimates for any $\ap\in (0,\bt\wedge(\sm_0-1))$.

\begin{definition} For $\sm\in(0,2)$ and an operator $\BI$, we
define the rescaled operator $\BI_{\mu,\ld}$ as in $(6.1)$. Then the
norm of scale $\sm$ is defined as
$$\|\BI^{(1)}-\BI^{(2)}\|_{\sm}=\sup_{(\mu,\ld)\in[1,\iy)\times(0,1)}\|\BI^{(1)}_{\mu,\ld}-\BI^{(2)}_{\mu,\ld}\|$$
where $\|\cdot\|$ is the norm defined in Definition
2.1.\end{definition}
\begin{remark} From $(6.2)$, we see that
$\|\BI^{(1)}-\BI^{(2)}\|_{\sm}\cong\sup_{\ld\in(0,1)}\|\BI^{(1)}_{1,\ld}-\BI^{(2)}_{1,\ld}\|$.
\end{remark}

The rescaled operator implies that if $u$ solves the equation $\BI
u-\pa_t u=f$ in $Q_{\ld}$, then the rescaled function
$w_{\mu}(x,t)=\mu u(\ld x,\ld^{\sm}t)$ solves an equation of the
same ellipticity type $\BI_{\mu,\ld}w_{\mu}(x,t)-\pa_t
w_{\mu}(x,t)=\ld^{\sm}\mu f(\ld x,\ld^{\sm}t)$ in $Q_1$.

The following theorem is the main result of this paper.

\begin{thm} Let $\sm\in(\sm_0,2)$ for $\sm_0\in(1,2)$ and let
$\BI^{(0)}$ be a fixed translation-invariant nonlocal operator in a
class $\fL\subset\fL_0(\sm)$ with scale $\sm$. Suppose that the
equation $\BI^{(0)}u-\pa_t u=0$ in $Q_{1+\e}$ has interior
$C^{1,\bt}$-estimates. Let $\BI^{(1)}$ and $\BI^{(2)}$ be two
nonlocal operators which are uniformly elliptic with respect to
$\fL_0(\sm)$ and assume that $\|\BI^{(0)}-\BI^{(k)}\|_{\sm}<\dt$ for
some $\dt>0$ small enough and $k=1,2$. If $\,u\in
L^{\iy}_T(L^1_{\om})$ solves the equations
$$\BI^{(1)}u-\pa_t u\ge f_1\,\text{ and }\,\BI^{(2)}u-\pa_t u\le
f_2\,\,\text{ in $Q_{1+\e}$}$$ for functions
$f_1,f_2\in\rB(Q_{1+\e})$, then $u\in C^{1,\ap}(Q_1)$ for any
$\ap\in(0,\bt\wedge(\sm_0-1))$, and moreover we have the estimate
\begin{equation}\begin{split}\|u\|_{C^{1,\ap}(Q_1)}&\lesssim
\|u\|_{C(Q_{1+\e})}+\|u\|_{L^{\iy}_T(L^1_{\om})}+(\sup_{Q_{1+\e}}|f_1|)\vee(\sup_{Q_{1+\e}}|f_2|).
\end{split}\end{equation}
\end{thm}

\pf We note that $u$ is continuous on $\overline{Q}_{1+\e}$. We
write $u=v+w$ where $v=u\mathbbm{1}_{Q_{1+\e}}$ and
$w=u\mathbbm{1}_{\BR^n_T\s Q_{1+\e}}$. Then by the uniform
ellipticity of $\BI^{(1)}$ we easily have that
\begin{equation*}\begin{split}\bM^+_{\fL_0}v-\pa_t v\ge-\|u\|_{L^{\iy}_T(L^1_{\om})}
-(\sup_{Q_{1+\e}}|f_1|)\vee(\sup_{Q_{1+\e}}|f_2|)\,\,\,\text{ in
$Q_1$. }
\end{split}\end{equation*} Similarly, we have that
\begin{equation*}\bM^-_{\fL_0}v-\pa_t
v\le\|u\|_{L^{\iy}_T(L^1_{\om})}+(\sup_{Q_{1+\e}}|f_1|)\vee(\sup_{Q_{1+\e}}|f_2|)\,\,\,\text{
in $Q_1$. }
\end{equation*} So we might use $v$ instead of $u$.

We now select some $\ld>0$ small enough so that
\begin{equation}\begin{split}&\ld^{\bt-\ap}+C_0 2^{\f{(1-\sm)(\bt-\ap_1)}{\sm}}\ld^{\ap_1-\ap}2^{\f{1+\ap_1}{\sm}}
<1,\\&\ld^{\bt-\ap}+2^{\f{(\sm-1)(1+\ap_1)}{\sm}}\ld^{\ap_1-\ap}(3+C_{\bt})<1,
\end{split}\end{equation}
where $\ap_1\in(\ap,\bt\wedge(\sm_0-1))$.
Take any $\vep>0$ with $\vep<\ld^{1+\bt}$. Then we choose $\dt=\dt(\vep)>0$
small enough as in Lemma 5.7. By scaling, without loss of
generality, we may assume that
$$\|u\|_{L^{\iy}_T(L^1_{\om})}+(\sup_{Q_{1+\e}}|f_1|)\vee(\sup_{Q_{1+\e}}|f_2|)<\dt\,\,\,\text{ and
}\,\,\,\sup_{\BR^n}|u|\le 1$$ and $u$ solves the equation in some
large cube $Q_R$.

By [C], it suffices to show that there are some $\ld\in(0,1)$ and a
sequence of linear functions
$$\el_k(x,t)=a_k+\la b_k,x\ra$$ such that
\begin{equation}\begin{split}
\sup_{Q_{\ld^k}}|u-\el_k|&\le\ld^{k(1+\ap)},\\
|a_{k+1}-a_k|&\le\ld^{k(1+\ap)},\\
\ld^k |b_{k+1}-b_k|&\le c_2\ld^{k(1+\ap)}.\end{split}\end{equation}
Set $\el_0=0$. Then we note that $|u_0|\le 1$ in $Q_1$ and
$|u_0(x,t)|\le (|x|^{\sm}+|t|)^{\f{1+\ap_1}{\sm}}$ for any
$(x,t)\in\BR^n_T\s Q_1$. We now continue the proof by the
mathematical induction. Assume that (6.5) holds for $k$-step. We
shall show that they are still working for $(k+1)$-step. We set
$$u_k(x,t)=\f{u(\ld^k x,\ld^{k\sm}t)-\el_k(\ld^k
x,\ld^{k\sm}t)}{\ld^{k(1+\ap)}}.$$

Since the class $\fL$ has scale $\sm$, $u_k$ solves equations of the
same ellipticity type as follows;
\begin{equation}\begin{split}
\BI_k^{(1)}u_k(x,t)&:=\BI^{(1)}_{\ld^{-k(1+\ap)},\ld^k}u_k(x,t)\ge\ld^{k(\sm-1-\ap)}f_1(\ld^k
x,\ld^{k\sm}t),\\
\BI_k^{(2)}u_k(x,t)&:=\BI^{(2)}_{\ld^{-k(1+\ap)},\ld^k}u_k(x,t)\le\ld^{k(\sm-1-\ap)}f_2(\ld^k
x,\ld^{k\sm}t).\end{split}\end{equation} We observe that the right
hand side (6.6) is getting smaller as $k$ increases. Thus we have
that
$$\|\BI_k^{(i)}-\BI_k^{(0)}\|\le\|\BI^{(i)}-\BI^{(0)}\|_{\sm}<\dt.$$
Let $\ap_1\in(\ap,\bt\wedge(\sm_0-1))$ be given. By the inductive
assumption, we see that $|u_k|\le 1$ in $Q_1$ and $|u_k(x,t)|\le
(|x|^{\sm}+|t|)^{\f{1+\ap_1}{\sm}}$ for any $(x,t)\in\BR^n_T\s Q_1$.
Then we shall construct functions $\el_{k+1}(x,t)$ and
$u_{k+1}(x,t)$ so that
$$|u_{k+1}(x,t)|\le (|x|^{\sm}+|t|)^{\f{1+\ap_1}{\sm}}\,\,\text{ for any $(x,t)\in\BR^n_T\s Q_1$.
}$$ Since $u$ is uniformly continuous on $\overline{Q}_{1+\e}$, we
may take some $R=R(\vep)>0$ (as in Lemma 5.7) so that $u$ admits a
modulus of continuity $\vr$ satisfying
$$|u(x,t)-u(y,s)|\le\vr((|x-y|^{\sm}+|t-s|)^{1/\sm})$$ for any $(x,t)\in
(Q_R\s Q_1)$ and $(y,s)\in \BR^n_T\s Q_1$. Then we apply Lemma 5.7
to the function $g$ which solves
\begin{equation*}\begin{cases} \BI_k^{(0)} g-\pa_t g=0 &\text{ in $Q_1$, }\\
g=u_k &\text{ in $\BR^n_T\s Q_1$ }\end{cases}\end{equation*} to
obtain that $\sup_{Q_1}|u_k-g|<\vep$. By the assumption, we note
that $\BI_k^{(0)}$ has interior $C^{1,\bt}$-estimates. Let $\hat
g(x,t)=\hat a+\la\hat b,x\ra$ be the linear part of $g$ at the
origin. Then we see that $|\hat a|<1+\vep$, because
$\sup_{Q_1}|g|\le 1+\vep$. By the $C^{1,\bt}$-estimates of $g$, we
have that
\begin{equation}|g(x,t)-\hat g(x,t)|\le
C_0(|x|^{\sm}+|t|)^{\f{1+\bt}{\sm}}\end{equation} for any $(x,t)\in
Q_{1/2}$. Since $\displaystyle\hat g\bigl(\f{\hat b}{2|\hat
b|},0\bigr)=\hat a+\f{1}{2}|\hat b|$, by (6.7) we have the upper
bound of $\hat b$ as follows;
\begin{equation*}\begin{split}|\hat b|&\le 2\bigl|\hat g\bigl(\f{\hat b}{2|\hat
b|},0\bigr)-g\bigl(\f{\hat b}{2|\hat b|},0\bigr)\bigr|
+2\bigl|g\bigl(\f{\hat b}{2|\hat b|},0\bigr)\bigr|+2|\hat a|\\
&=C_0 2^{-\bt}+4(1+\vep):=C_{\bt}.
\end{split}\end{equation*}
Then we can derive the estimates as follows;
\begin{equation}\begin{split}
|u_k(x,t)-\hat g(x,t)|\le\begin{cases}\vep+C_0(|x|^{\sm}+|t|)^{\f{1+\bt}{\sm}},&(x,t)\in Q_{1/2},\\
\vep+2+C_\bt,&(x,t)\in Q_1\s Q_{1/2},\\
\vep+(2+C_\bt)(|x|^{\sm}+|t|)^{\f{1+\ap_1}{\sm}},&(x,t)\in\BR^n_T\s
Q_1.\end{cases}
\end{split}\end{equation}
We now set
\begin{equation}\begin{split}\el_{k+1}(x,t)&=\el_k(x,t)+\ld^{k(1+\ap)}\hat
g(\ld^{-k}x,\ld^{-k\sm}t),\\
u_{k+1}(x,t)&:=\f{u(\ld^{k+1} x,\ld^{(k+1)\sm}t)-\el_{k+1}(\ld^{k+1}
x,\ld^{(k+1)\sm}t)}{\ld^{(k+1)(1+\ap)}}\\
&\,\,=\f{u_k(\ld x,\ld^{\sm}t)-\hat g(\ld
x,\ld^{\sm}t)}{\ld^{1+\ap}}.
\end{split}\end{equation}
Then it follows from (6.5) and (6.9) that
\begin{equation}\begin{split}
&|u_{k+1}(x,t)|\\&\le\begin{cases}\ld^{\bt-\ap}+C_0
2^{\f{1-\sm}{\sm}(\bt-\ap_1)}\ld^{\ap_1-\ap}(|x|^{\sm}+|t|)^{\f{1+\ap_1}{\sm}},
&(x,t)\in Q_{\ld^{-1}/2},\\
2^{\f{\sm-1}{\sm}(1+\ap_1)}\ld^{\ap_1-\ap}(3+C_{\bt})(|x|^{\sm}+|t|)^{\f{1+\ap_1}{\sm}},
&(x,t)\in Q_{\ld^{-1}}\s Q_{\ld^{-1}/2},\\
\ld^{\bt-\ap}+\ld^{\ap_1-\ap}(2+C_{\bt})(|x|^{\sm}+|t|)^{\f{1+\ap_1}{\sm}},&(x,t)\in\BR^n_T\s
Q_{\ld^{-1}},\end{cases}
\end{split}\end{equation}
and moreover $|u_{k+1}|\le 1$ on $Q_1$ and $|u_{k+1}(x,t)|\le
(|x|^{\sm}+|t|)^{\f{1+\ap_1}{\sm}}$ for any $(x,t)\in\BR^n_T\s Q_1$.
Finally, it follows from (6.10) that
$$u(x,t)-\el_{k+1}(x,t)=\ld^{(k+1)(1+\ap)}u_{k+1}(\ld^{-(k+1)}x,\ld^{-(k+1)\sm}t)$$
for any $(x,t)\in Q_{\ld^{k+1}}$, and thus we conclude that
\begin{equation*}\sup_{Q_{\ld^{k+1}}}|u-\el_{k+1}|\le\ld^{(k+1)(1+\ap)}\sup_{Q_1}|u_{k+1}|\le\ld^{(k+1)(1+\ap)}.
\end{equation*} Hence we complete the proof. \qed

\section{Cordes-Nirenberg type estimates and Applications}

In this section, we furnish various applications of the previous results.
Our proofs on the parabolic case are based on that
\cite{CS1} of the elliptic case and the results \cite{KL} of the
parabolic case, but the main difference between them is in that we
extend the solution space $\rB(\BR^n)$ on the elliptic one to the
more flexible space $L^{\iy}_T(L^1_{\om})$ on the parabolic one and
use the more wider class of kernels involving variables
$(x,t)\in\BR^n_T$, and moreover the parabolic case requires more
careful consideration due to the time shift. Contrary to the
elliptic case, we had better mention on the difficulty in the
parabolic case; for instance, "time shift".

$$\text{7.1. a parabolic version of the integral
Cordes-Nirenberg type estimates}$$

When the equation is linear and close to an operator in $\fL_1$ in
an appropriate way, we shall obtain the regularity results of its
viscosity solutions. This is a parabolic version of the integral
Cordes-Nirenberg type estimates.

\begin{thm} For $\sm\in(1,2)$, let $u\in L^{\iy}_T(L^1_{\om})$ be a
viscosity solution of the equation $$\BI u-\pa_t
u:=\int_{\BR^n}\mu_{\cdot}(u,\cdot,y)\f{(2-\sm)a(\,\cdot\,,y,\,\cdot\,)}{|y|^{n+\sm}}\,dy-\pa_t
u=f\,\text{ in $Q_{1+\e}$,}$$ where $f\in\rB(Q_{1+\e})$. Suppose
that there is some $\dt>0$ small enough such that
\begin{equation}\sup_{y\in\BR^n}\sup_{(x,t)\in
Q_3}|a(x,y,t)-a_0(x,y,t)|<\dt,\end{equation} where $a_0$ is a
bounded function so that $(2-\sm)a_0(x,y,t)/|y|^{n+\sm}$ satisfies
$(1.2)$.

If $u\in L^{\iy}_T(L^1_{\om})$ is a viscosity solution of the
equation
$$\int_{\BR^n}\mu_{\cdot}(u,\cdot,y)\f{(2-\sm)a_0(\,\cdot\,,y,\,\cdot\,)}{|y|^{n+\sm}}\,dy-\pa_t
u=0\,\,\text{ in $Q_{1+\e}$,}$$ then there is some $\bt\in(0,1)$ so
that $u\in C^{1,\ap}(Q_1)$ for any $\ap\in(0,\bt\wedge(\sm-1))$ and we have the estimate
$$\|u\|_{C^{1,\ap}(Q_1)}\lesssim
\|u\|_{C(Q_{1+\e})}+\|u\|_{L^{\iy}_T(L^1_{\om})}+\sup_{Q_{1+\e}}|f|.$$
\end{thm}

\pf Without loss of generality, we may assume that
$\|u\|_{C(Q_{1+\e})}\vee\|u\|_{L^{\iy}_T(L^1_{\om})}\lesssim 1$ by dividing $u$ by
$\|u\|_{C(Q_{1+\e})}+\|u\|_{L^{\iy}_T(L^1_{\om})}+\sup_{Q_{1+\e}}|f|$.
We apply Theorem 6.3. In this case, $\BI^{(0)}$ is given by
$$\BI^{(0)}u(x,t)=\int_{\BR^n}\mu_t(u,x,y)\f{(2-\sm)a_0(x,y,t)}{|y|^{n+\sm}}\,dy.$$
By the assumption, the operator $\BI^{(0)}$ is translation-invariant
and belongs to $\fL_1$, which is a scale invariant class. By Theorem
3.6, the viscosity solution $u\in L^{\iy}_T(L^1_{\om})$ of
$\BI^{(0)}u-\pa_t u=0$ in $Q_{1+\e}$ has interior
$C^{1,\bt}(Q_1)$-estimates for some $\bt\in(0,1)$. In this case, it
is easy to see that since the equation is linear and the
coefficients do not depend on $(x,t)$, the derivatives $Du$ of the
solution $u$ of the equation solve the same equation, so that the
solutions are actually $C_x^{2,\bt}(Q_1)$. So, moreover, the
solutions $u$ of $\BI^{(0)}u-\pa_t u=0$ has interior
$C_x^{1,1}(Q_1)$-estimates.

Now we estimate
$\|\BI-\BI^{(0)}\|_{\sm}:=\sup_{\ld\in(0,1)}\|\BI_{1,\ld}-\BI^{(0)}_{1,\ld}\|$.
We take any $\ld\in(0,1)$ and set $w_1=u(\ld\,\cdot,\ld^{\sm}\cdot)$
in $Q_{1+\e}$. Take any $(y,s)\in Q_{(1+\e)\ld}$. Then we compute
$$(\BI_{1,\ld}-\BI^{(0)}_{1,\ld})u(y,s)=\int_{\BR^n}\mu_s(u,y,z)\f{(2-\sm)\bigl(a(y,\ld z,s)
-a_0(y,\ld z,s)\bigr)}{|z|^{n+\sm}}\,dz$$ by applying (7.1). By
Definition 3.1, we have that
\begin{equation}\begin{split}\|\BI_{1,\ld}-\BI^{(0)}_{1,\ld}\|&=\sup_{(y,s)\in
Q_{(1+\e)\ld}}\sup_{u\in\cF^{M}_{y,s}}
\f{|\BI_{1,\ld}u(y,s)-\BI^{(0)}_{1,\ld}u(y,s)|}{1+\|u\|_{L^{\iy}_T(L^1_{\om})}+\|u\|_{C^{1,1}_x(Q_1(y,s))}}\\
&\le\sup_{(y,s)\in
Q_{(1+\e)\ld}}\sup_{u\in\cF^{M}_{y,s}}\f{\ds\int_{\BR^n}|\mu_s(u,y,z)|\f{(2-\sm)\,\dt}{|z|^{n+\sm}}\,dz}
{1+\|u\|_{L^{\iy}_T(L^1_{\om})}+\|u\|_{C^{1,1}_x(Q_1(y,s))}}
\end{split}\end{equation} We take any function $u\in\cF^{M}_{y,s}$;
that is,
 $\|u\|_{L^{\iy}_T(L^1_{\om})}+\|u\|_{C^{1,1}_x(Q_1(y,s))}\le M$
for $M>0$ and $u\in\fF\cap C^2_x(y,s)$. Since $y\pm z\in B_1(y)$ if
$\pm z\in B_1$, we obtain that
\begin{equation*}\begin{split}\int_{B_1}|\mu_s(u,y,z)|\f{2-\sm}{|z|^{n+\sm}}\,dz
&\le\int_{B_1}\|u\|_{\rC^{1,1}_x(Q_1(y,s))}|z|^2\f{2-\sm}{|z|^{n+\sm}}\,dz\lesssim
\|u\|_{C^{1,1}(Q_1(y,s))}.\end{split}\end{equation*} Since
$(y,s)\in Q_{\ld(1+\e)}$ and $|y\pm z|\ge
|z|-|y|\ge(1-\ld-\ld\e)|z|$ for $y\in B_{(1+\e)\ld}$ and
$z\in\BR^n\s B_1$, we have that
\begin{equation*}\begin{split}&\int_{\BR^n\s
B_1}|\mu_s(u,y,z)|\f{2-\sm}{|z|^{n+\sm}}\,dz\\
&\qquad\qquad\le\int_{\BR^n\s B_1}\bigl(|u(y+z,s)|+|u(y-z,s)|+2|u(y,s)|\bigr)\f{2-\sm}{|z|^{n+\sm}}\,dz\\
&\qquad\qquad=\int_{|z|\ge
1-\ld-\ld\e}|u(z,s)|\bigr(\f{2-\sm}{|y+z|^{n+\sm}}+\f{2-\sm}{|y-z|^{n+\sm}}\bigr)\,dz+C\\
&\qquad\qquad\lesssim\int_{|z|\ge
1-\ld-\ld\e}|u(z,s)|\f{2-\sm}{|z|^{n+\sm}}\,dz+1\lesssim
\|u\|_{L^{\iy}_T(L^1_{\om})}+1\lesssim 1.
\end{split}\end{equation*}
By (7.2), we conclude that $\|\BI_{1,\ld}-\BI^{(0)}_{1,\ld}\|\le
C\dt$ for any $\ld\in(0,1)$, and thus we have that
$\|\BI-\BI^{(0)}\|_{\sm}\le C\dt$. If we choose $\e$ small enough,
we can apply Theorem 6.3 and conclude that the equation $\BI u-\pa_t
u=f$ has interior $C^{1,\ap}$-estimates for any
$\ap\in(0,\bt\wedge(\sm-1))$.\qed

$$\text{7.2. Nonlinear equations}$$

From Theorem 6.3 and Theorem 3.6, we can easily derive the following
result.

\begin{thm} Let $\sm\in[\sm_0,2)$ for $\sm_0\in(1,2)$ and
let $\fL_1(\sm)$ be the class satisfying $(1.3)$. Suppose that $u\in
L^{\iy}_T(L^1_{\om})$ is a viscosity solution of the equation
$$\BI^{(0)}u-\pa_t u=0\,\,\text{ in $Q_{1+\e}$ }$$ and
$\|\BI-\BI^{(0)}\|_{\sm}<\dt$ for some small $\dt>0$, where
$\BI^{(0)}$ is a translation-invariant nonlocal operator which is
uniformly elliptic with respect to $\fL_1(\sm)$ and $\BI$ is an
operator which is uniformly elliptic with respect to $\fL_0$. If
$u\in L^{\iy}_T(L^1_{\om})$ is a viscosity solution of the equation
$$\BI u-\pa_t u=f\,\text{ in $Q_{1+\e}$}$$ where $f\in\rB(Q_{1+\e})$, then $u\in C^{1,\ap}(Q_1)$
for some small $\ap>0$ and we have the estimate $$\|u\|_{C^{1,\ap}(Q_1)}\lesssim
\|u\|_{C(Q_{1+\e})}+\|u\|_{L^{\iy}_T(L^1_{\om})}+\sup_{Q_{1+\e}}|f|.$$
\end{thm}

\begin{remark} We consider the following operator $\BI$
given by $$\BI
u:=\inf_{\ap}\sup_{\bt}\int_{\BR^n}\mu_{\cdot}(u,\cdot,y)\f{(2-\sm)\bigl(a_0(\,\cdot\,,y,\,\cdot\,)
+a_{\alpha\beta}(\,\cdot\,,y,\,\cdot\,)\bigr)} {|y|^{n+\sm}}\,dy$$
where $a_0$ and $a_{\alpha\beta}$ are functions satisfying
\begin{equation*}\begin{split}\ld\le a_0(x,y,t)\le&\Ld,\,\,\sup_{(x,t)\in\BR^n_T}|\n_y a_0(x,y,t)|\le\f{C}{|y|}\,\text{
for any $y\in\BR^n\s\{0\}$,
}\\\sup_{\ap,\bt}&\sup_{(x,t)\in\BR^n_T}|a_{\alpha\beta}(x,y,t)|<\dt\,\,\,\text{
for some small $\dt>0$.}\end{split}\end{equation*} Then we see that
this is a nonlinear operator which exemplifies Theorem 7.2.
\end{remark}

\begin{thm} Let $\sm\in[\sm_0,2)$ for $\sm_0\in(1,2)$ and
let $\BI u$ be given by
$$\BI u=\inf_{\ap}\sup_{\bt}\int_{\BR^n}\mu_{\cdot}(u,\cdot,y)\f{(2-\sm)a_{\alpha\beta}(\,\cdot\,,y,\,\cdot\,)}{|y|^{n+\sm}}\,dy$$
where $\ld<a_{\alpha\beta}(x,y,t)<\Ld$, $|\n_y a_{\alpha\beta}(x,y,t)|\le
C_2/|y|$ and
$$\sup_{\ap,\bt}|a_{\alpha\beta}(x_1,y,t_1)-a_{\alpha\beta}(x_2,y,t_2)|=o(1)
\,\text{ as $\dd((x_1,t_1),(x_2,t_2))\to 0$ }$$ with the parabolic
distance $\dd$. If $u\in L^{\iy}_T(L^1_{\om})$ is a viscosity
solution of the equation
$$\BI u-\pa_t u=f\,\text{ in $Q_{1+\e}$ }$$ where $f\in\rB(Q_{1+\e})$,
then there is a small $\ap>0$ and we have the estimate
$$\|u\|_{C^{1,\ap}(Q_1)}\lesssim
\|u\|_{C(Q_{1+\e})}+\|u\|_{L^{\iy}_T(L^1_{\om})}+\sup_{Q_{1+\e}}|f|.$$
\end{thm}

\pf For each $(x_0,t_0)\in Q_1$, we can find a ball
$Q_r(x_0,t_0)\subset Q_{1+\e}$ ($r>0$ is independent of $x_0$ and
$t_0$) so that
$$\sup_{(x,t)\in Q_r(x_0,t_0)}|a_{\alpha\beta}(x,y,t)-a_{\alpha\beta}(x_0,y,t_0)|<\dt$$ for
some small $\dt>0$. This implies that
$\|\BI-\BI_{(x_0,t_0)}\|_{\sm}<C\dt$ on $Q_r(x_0,t_0)$ as in the
proof of Theorem 7.2, where
\begin{equation*}\begin{split}\BI_{(x_0,t_0)}u(x,t)&:=\inf_{\ap}\sup_{\bt}\int_{\BR^n}\mu_t(u,x,y)
\f{(2-\sm)a_{\alpha\beta}(x_0,y,t_0)}{|y|^{n+\sm}}\,dy\\
&=\BI\bigl(\btau_{x-x_0}^{t-t_0}u\bigr)(x_0,t_0).\end{split}\end{equation*}
We now apply Theorem 7.2 with $\BI^{(0)}=\BI_{(x_0,t_0)}$ scaled in
$Q_r(x_0,t_0)$. Let $N$ be the minimal number of such open balls
$Q_r(x_0,t_0)$ covering $\overline Q_{1+\e}$. Then we have that
$$\|u\|_{C^{1,\ap}(Q_{1+\e})}\le
C_N\bigl(\,\|u\|_{C(Q_{1+\e})}+\|u\|_{L^{\iy}_T(L^1_{\om})}+\sup_{Q_{1+\e}}|f|\bigr).$$
Hence we complete the proof. \qed

$$\text{7.3. Nonlinear equations with
non-differentiable kernels }$$

We note that Theorem 6.3 makes it possible to obtain certain results
even in the translation-invariant case. It was crucial in Theorem
3.6 that every kernel must be differentiable away from the origin.
This condition can be weakened in the following way. We establish
$C^{1,\ap}$-estimates for nonlocal equations that are uniformly
elliptic with respect to the class $\fL$ consisting of operators
with kernels $K\in\cK$ given by
$$K(x,y,t)=(2-\sm)\f{a_1(x,y,t)+a_2(x,y,t)}{|y|^{n+\sm}}$$ where $\ld\le
a_1\le\Ld$,
$$\sup_{y\in\BR^n}\sup_{(x,t)\in\BR^n_T}|a_2(x,y,t)|<\dt,\,\sup_{(x,t)\in\BR^n_T}|\n_y
a_1(x,y,t)|\le\f{c_1}{|y|}\,\,\text{ for any $y\in\BR^n\s\{0\}$.} $$

\begin{thm} Let $\sm\in[\sm_0,2)$ for $\sm_0\in(1,2)$ and let
$\dt>0$ be a small enough number $($depending only on $\ld,\Ld,c_1$
and the dimension $n$, but not on $\sm$$)$ as in the above. If $u\in
L^{\iy}_T(L^1_{\om})$ is a viscosity solution of the nonlocal equation
$$\BI u-\pa_t u:=\inf_{\ap}\sup_{\bt}\int_{\BR^n}\mu_{\cdot}(u,\cdot,y)K_{\ap\bt}(y,\cdot)\,dy-\pa_t
u=f\,\text{ in $Q_{1+\e}$ }$$ for $f\in\rB(Q_{1+\e})$ and
$\{K_{\ap\bt}\}\subset\cK_0$, then there is some $\ap>0$ such that
$$\|u\|_{C^{1,\ap}(Q_1)}\lesssim
\|u\|_{C(Q_{1+\e})}+\|u\|_{L^{\iy}_T(L^1_{\om})}
+\sup_{Q_{1+\e}}|f|.$$
\end{thm}

\pf Let $L\in\fL$ be an operator with kernel $K$. We write
$K=K_1+K_2$ where $K_1=(2-\sm)a_1(x,y,t)/|y|^{n+\sm}$ and set
$L=L^1+L^2$ where $L^1$ and $L^2$ are operators with kernels $K_1$
and $K_2$, respectively. Then we see that $\|L-L^1\|_{\sm}<c\dt$. If
we set $\BI^{(0)}u=\inf_{\ap}\sup_{\bt}L^1_{\ap\bt}u$, then we have
that $\|\BI-\BI^{(0)}\|_{\sm}<c\dt$, and hence we can apply Theorem
6.3 to complete the proof. \qed

\rk This theorem works for a class which is still much smaller than
$\fL_0$. It would be very interesting to determine whether the class
$\fL_0$ has interior $C^{1,\ap}$-estimates or not. This problem is
still left open even for elliptic cases as mentioned in \cite{CS1}.
Also it would be interesting to answer this problem on the parabolic
case.

$$\text{7.4. Nonlinear equations near the fractional Laplacian}$$

We obtain another result in translation-invariant case by applying
Theorem 6.3. In fact, we obtain $C^{2,\ap}$-estimates for nonlinear
translation-invariant nonlocal parabolic equations which are
sufficiently close to the parabolic fractional Laplacian and their
ellipticity constants are sufficiently close to each other. This is
to improve Theorem 3.6 under these conditions.

\begin{thm} Let $\sm\in[\sm_0,2)$ for $\sm_0\in(1,2)$. Then
there are some $\dt>0$ and $\rho_0>0$ so that if
$-\dt<\ld-1<\Ld-1<\dt$, $\BI$ is a nonlocal translation-invariant
uniformly elliptic operator with respect to $\fL_*$ and $u\in
L^{\iy}_T(L^1_{\om})$ is a viscosity solution of the equation $\BI u-\pa_t u=0$ in
$Q_{1+\e}$, then $u\in C_x^{2,\ap}(Q_1)$ for a constant $\ap\in(0,1)$
$($depending only on $n$ and $\sm_0$$)$ and we have the estimate
$$\|u\|_{C_x^{2,\ap}(Q_1)}\lesssim
\|u\|_{L^{\iy}_T(L^1_{\om})}+|\BI 0|$$ where we denote by $\BI 0$ the value we obtain when we apply $\BI$ to the constant function that is equal to zero.\end{thm}

\pf From Theorem 3.6, we see that $u\in C^{1,\ap}(Q_1)$. Thus the
function $u$ is differentiable in $x$ on $Q_1$. Let $w=e\cdot\n u$
be a directional derivative for $e\in S^{n-1}$. We write $w=w_1+w_2$
where $w_1=w\mathbbm{1}_{Q_{1+\e}}$. Then by using the uniform
ellipticity with respect to $\fL_*$ we easily see that $w_1$ solves
$$\BM^+_{\fL_*}w_1-\pa_t w_1\ge-\|u\|_{L^{\iy}_T(L^1_{\om})}-|\BI 0|,\,\,\BM^-_{\fL_*}w_1-\pa_t w_1\le
\|u\|_{L^{\iy}_T(L^1_{\om})}+|\BI 0|\text{ in $Q_{1+2\e/3}$. }$$ We
now apply Theorem 6.3 instead of Theorem 3.4. Since
$1-\dt<\ld<\Ld<1+\dt$, as in (7.2) we easily obtain that
$$\|\BM^+_{\fL_*}+(-\Delta)^{\sm/2}\|_{\sm}<c\dt\,\,\text{ and }
\,\,\|\BM^-_{\fL_*}+(-\Delta)^{\sm/2}\|_{\sm}<c\dt.$$ Thus Theorem
6.3 tells us that $w=e\cdot\n u$ is in $C_x^{1,\ap}(Q_1)$. From
(3.6) and the local equivalence between $W^{1,\iy}$ and Lipschitz
continuity, we see that $\sup_{Q_{1+\e}}|\n u|\le
\|u\|_{C^{0,1}_x(Q_{1+\e})}\lesssim\|u\|_{L^{\iy}_T(L^1_{\om})}$.
Also by (3.7) and integration by parts, we have that $\|\n u
\|_{L^{\iy}_T(L^1_{\om})}\lesssim\|u \|_{L^{\iy}_T(L^1_{\om})}$.
Moreover, if we take $e=\n u/|\n u|$ in the above, then we have that
\begin{equation*}\begin{split}\|\n u\|_{C_x^{1,\ap}(Q_1)}&\le\|\n u\|_{C^{1,\ap}(Q_1)}\\
&\lesssim\|\n u\|_{C(Q_{1+\e})}+\|\n u\|_{L^{\iy}_T(L^1_{\om})}+|\BI 0|\\
&\lesssim\|u \|_{L^{\iy}_T(L^1_{\om})}+|\BI 0|.
\end{split}\end{equation*}
This implies that $\|u\|_{C_x^{2,\ap}(Q_1)}\lesssim\|u\|_{L^{\iy}_T(L^1_{\om})}
+|\BI 0|.$ \qed

When we consider the concave equation $\bI_0 u-\pa_t u=0$ in $Q_2$ where $\bI_0$ is defined in $\fL=\fL_2$, its viscosity solution $u\in L^{\iy}_T(L^1_{\om})$ admits a parabolic $C^{2,\ap}(Q_1)$-estimate. For this estimate, we need a lemma which can be shown as in the proof of Theorem 3.6.

\begin{lemma} Let $\sm\in[\sm_0,2)$ for some $\sm_0\in(1,2)$.
If $u\in L^{\iy}_T(L^1_{\om})$ is a viscosity solution of the parabolic equation
$$\BI u-\pa_t u=c_0\,\,\text{ in $Q_{1+\e}$ }$$ where $\BI$ is a nonlocal, translation-invariant and uniformly elliptic with respect to
$\fL_*$ and $c_0\in\BR$ is a constant,
then there is some $\ap\in(0,1)$ such that
$$\|u\|_{C^{1,\ap}(Q_1)}\lesssim
\|u\|_{L^{\iy}_T(L^1_{\om})}.$$\end{lemma}

\pf It can be proved in the almost same way as Theorem 3.6. \qed

\begin{thm} Let $\sm\in[\sm_0,2)$ for $\sm_0\in(1,2)$. Then there are some $\dt>0$ and $\rho_0>0$ so that if $-\dt<\ld-1<\Ld-1<\dt$ and $u\in
L^{\iy}_T(L^1_{\om})$ is a viscosity solution of the parabolic equation 
\begin{equation}\BI u-\pa_t u=0\,\,\text{ in
$Q_2$ }\end{equation} where $\BI$ is a nonlocal, translation-invariant and uniformly elliptic with respect to
$\fL_*$, then $u\in C_x^{2,\ap}(Q_1)\cap C_t^{1,\f{2+\ap-\sm}{\sm}}(Q_1)$ for a constant $\ap\in(0,1)$
$($depending only on $n$ and $\sm_0$$)$ and we have the estimate
$$\|u\|_{C_x^{2,\ap}(Q_1)}+\|u\|_{C_t^{1,\f{2+\ap-\sm}{\sm}}(Q_1)}\lesssim
\|u\|_{L^{\iy}_T(L^1_{\om})}.$$
\end{thm}

\pf Without loss of generality, we may assume that $\|u\|_{L^{\iy}_T(L^1_{\om})}=1$ by dividing $u$ by  $\|u\|_{L^{\iy}_T(L^1_{\om})}$.
As in Theorem 7.6, we obtain the estimate
\begin{equation}\|u\|_{C_x^{2,\ap}(Q_1)}\lesssim
\|u\|_{L^{\iy}_T(L^1_{\om})}=1
\end{equation} for some $\ap\in(0,1)$. Since
$(2+\ap)/\sm>1$ for such $\ap>0$, we see that
$$\f{2+\ap-\sm}{\sm}+1=\f{2+\ap}{\sm}$$ and $0<\ap<2+\ap-\sm<1$. Thus it suffices to show that $u$ is differentiable in $t$ on $Q_1$ and admits the estimate
\begin{equation}\|u\|_{C_t^{1,\f{2+\ap-\sm}{\sm}}(Q_1)}\lesssim
\|u\|_{L^{\iy}_T(L^1_{\om})}=1.\end{equation}
In order to get the estimate (7.4), we proceed it as in the proof of Theorem 3.6.
For $(x_0,t_0)\in Q_1$ and $c_0\in\BR$, we consider the function
$$w(x,t)=\f{\btau^{t_0}_{x_0} u(rx,r^{\sm}t)-u(x_0,t_0)-r\n u(x_0,t_0)\cdot x-c_0 r^{\sm}t-\f{1}{2}x^T\cdot D^2 u(x_0,t_0)\cdot x}{r^{2+\ap}}$$
for any sufficiently small $r>0$. Here $c_0$ is a constant to be chosen later.
If we set $v(x,t)=\btau^{t_0}_{x_0} u(x,t)-c_0 t$, then it satisfies the equation
\begin{equation}\BI v-\pa_t v=c_0\,\,\text{ in $Q_1$ }
\end{equation} and we have that
\begin{equation}\begin{split}w(x,t)&=\f{v(rx,0)-v(0,0)-r\n v(0,0)\cdot x-\f{1}{2}x^T\cdot D^2 v(0,0)\cdot x}{r^{2+\ap}}\\
&\qquad+\f{v(rx,r^{\sm}t)-v(rx,0)}{r^{2+\ap}}:=w_1(x,t)+w_2(x,t).
\end{split}\end{equation}
By Lemma 7.7, we see that $v\in C^{1,\ap}(Q_1)$ for some $\ap\in(0,1)$. So the directional derivative $e\cdot\n v$ for $e\in S^{n-1}$ satisfies the equation (7.3).
Thus, as in the proof of Theorem 7.6, we have the estimate
\begin{equation}\|v\|_{C_x^{2,\ap}(Q_1)}\lesssim\|v\|_{L^{\iy}_T(L^1_{\om})}\lesssim\|u\|_{L^{\iy}_T(L^1_{\om})},\end{equation} and so 
\begin{equation}\sup_{Q_1}|w_1|\lesssim\|u\|_{L^{\iy}_T(L^1_{\om})}.
\end{equation}

We observe that $w_2$ solves the equation (7.6).
Considering the function $\phi(y)=|y|^2\mathbbm{1}_{B_{1+\e}}(y)+(1+\e)^2\mathbbm{1}_{\BR^n\s B_{1+\e}}$, as in the proof of Theorem 3.6 we have that $\BI\phi\le 6\Ld\om_n\e^{-\sm}$ on $B_1$. 
Put $M_*=\sup_{B_1\times(-1,t_*)}w_2$. Then without loss of generality we may assume that $M_*\ge 0$; otherwise, we could use $-w$ instead of $w$. If we take some $c_0<M_*$ and $M_*=w_2(x_*,t_*)\ge  0$ for some $x_*\in B_1$ and $t_*\in(-1,t_1)$ where $t_1=-1+\f{1}{12\Ld\om_n\e^{-\sm}}$, then it is not difficult to check that the functions
\begin{equation*}\begin{split}
\phi_3(x,t)&=(M_*-c_0)(t+1)+M_*\f{\|v\|_{C_x^{2,\ap}(Q_1)}\,\phi(x)}{12\Ld\om_n\e^{-\sm}(1+\e)^2},\\
\phi_4(x,t)&=6\Ld\om_n\e^{-\sm}(M_*-c_0)(t+1)+M_*\f{\phi(x)-\phi(x_*)}{6\Ld\om_n\e^{-\sm}}
+\|v\|_{C^{2,\ap}_x(Q_1)},
\end{split}\end{equation*} are supersolutions of the equation (7.6), provided that $\e$ could be chosen so small that $6\Ld\om_n\e^{-\sm}>M_*\vee 1$. As in the proof of Theorem 3.6, we can derive that $w_2\le M_*+\|v\|_{C_x^{2,\ap}(Q_1)}$ and $M_*\le 4\|v\|_{C_x^{2,\ap}(Q_1)}$, and thus by (7.8) we get that $w_2\lesssim\|u\|_{L^{\iy}_T(L^1_{\om})}$ on $Q_1$. Similarly, we can obtain the estimate $w_2\gtrsim-\|u\|_{L^{\iy}_T(L^1_{\om})}$ on $Q_1$ by constructing subsolutions corresponding to $\phi_3$ and $\phi_4$. This implies that 
\begin{equation}\sup_{Q_1}|w_2|\lesssim\|u\|_{L^{\iy}_T(L^1_{\om})}.
\end{equation}
By (7.9) and (7.10), we thus conclude that $\sup_{Q_1}|w|\lesssim\|u\|_{L^{\iy}_T(L^1_{\om})}$. If we set $s=r^{\sm}$, then this estimate gives that
\begin{equation*}\begin{split}
|u(x_0,t_0+s t)-u(x_0,t_0)-c_0 s t|=\cO(s^{\f{2+\ap}{\sm}})=\fo(s)
\end{split}\end{equation*}
for any sufficiently small $s>0$. Thus $u$ is differentiable in time at $(x_0,t_0)$, and moreover $c_0=\pa_t u(x_0,t_0)$. Hence we have that
\begin{equation*}|w(0,t)|=\biggl|\f{u(0,t_0+st)-u(0,t_0)-\pa_t u(0,t_0) st}{s^{\f{2+\ap}{\sm}}}\biggr|\lesssim\|u\|_{L^{\iy}_T(L^1_{\om})}
\end{equation*} for any sufficiently small $s>0$. This and the mean value theorem lead us to get the inequality (7.5). Therefore we complete the proof. \qed

From Remark 2.1 and Theorem 7.8, we can easily obtain the following corollary.

\begin{cor} Let $\sm\in[\sm_0,2)$ for $\sm_0\in(1,2)$. Then there are some $\dt>0$ and $\rho_0>0$ so that if $-\dt<\ld-1<\Ld-1<\dt$ and $u\in
L^{\iy}_T(L^1_{\om})$ is a viscosity solution of the parabolic equation 
\begin{equation}\bI_0 u-\pa_t u=0\,\,\text{ in
$Q_2$ }\end{equation} where $\bI_0$ is defined in $\fL=\fL_2(\sm)$, then $u\in C^{2,\ap}(Q_1)$ for a constant $\ap\in(0,1)$
$($depending only on $n$ and $\sm_0$$)$ and we have the estimate
$$\|u\|_{C^{2,\ap}(Q_1)}\lesssim\|u\|_{L^{\iy}_T(L^1_{\om})}.$$
\end{cor}

\end{document}